\documentclass[preprint]{elsarticle}
\usepackage{xcolor}
\usepackage{amsfonts}
\usepackage{amsmath}
\usepackage{amssymb}
\usepackage[title,titletoc]{appendix}
\usepackage{enumitem}
\usepackage[top=1.2in,bottom=1.2in,left=1in, right=1in]{geometry}
\usepackage{graphics}
\usepackage{lineno}
\usepackage{pgfplots}
\usepackage{tikz}
\usepackage{todonotes}
\usepackage{bm}
\usepackage{commath}
\usepackage{amsthm}
\usepackage{subfig}

\usepackage[pagebackref=false,bookmarks=false,hidelinks=true]{hyperref} 
\hypersetup{pdfborder = 0 0 0, colorlinks=false,
 linkcolor=black,citecolor=black, filecolor=black, urlcolor=black, }

\hypersetup{
  bookmarksnumbered=true,
  bookmarksopen=false,
  hypertexnames=false,      
  breaklinks=true,          
  unicode=false,
  pdffitwindow=true,        
  pdfnewwindow=true,        
  colorlinks=true,         
  linkcolor=dblue,
  anchorcolor=red,
  citecolor=dorange,
  filecolor=magenta,
  urlcolor=dblue,
  pdfstartview = FitH,
  pdfkeywords = {},
  pdfcreator = {LaTeX with hyperref package}
}

\usepackage{pgfplots,tikzscale}
\usetikzlibrary{plotmarks}
\usetikzlibrary{arrows}
\tikzset{mark size=1/2}

\usepackage{cleveref}
\crefname{equation}{}{}
\Crefname{equation}{}{}

\newtheorem{lemma}{Lemma}
\newtheorem{corollary}{Corollary}
\newtheorem{remark}{Remark}

\newcommand{\bigO}{\mathcal{O}}

\newcommand{\DD}{{\mathcal{D}}}

\newcommand{\ff}{{\bm{f}}}
\newcommand{\boldg}{{\bm{g}}}
\newcommand{\FF}{\bm{F}}

\newcommand{\kk}{{\bm{k}}}

\newcommand{\NN}{{\mathcal{N}}}

\newcommand{\pderiv}[2]{\frac{\partial #1}{\partial #2}}

\renewcommand{\Re}{\operatorname{\mathit{Re}}} % Reynolds

\newcommand{\rr}{{\bm{r}}}

\newcommand{\ssigma}{{\bm\sigma}}
\newcommand{\mmu}{{\bm\mu}}

\newcommand{\uu}{{\bm{u}}}

\newcommand{\zzero}{{\bm{0}}}
\newcommand{\vv}{{\bm{v}}}

\newcommand{\xx}{{\bm{x}}}
\newcommand{\xxi}{{\bm{\xi}}}
\newcommand{\yy}{{\bm{y}}}
\newcommand{\ZZ}{{\mathbb{Z}}}

\newcommand{\pfmm}{{p_{\textrm{FMM}}}}

\newcommand{\Amat}{{\bm A}}

\renewcommand{\v}[1]{\bm{#1}}
\newcommand{\di}{\partial_i}
\renewcommand{\dj}{\partial_j}
\newcommand{\dk}{\partial_k}

\newcommand{\nhat}{{\hat{n}}}
\newcommand{\that}{{\hat{t}}}
\newcommand{\real}{\operatorname{Re}}
\newcommand{\imag}{\operatorname{Im}}
\newcommand{\domain}{\Omega}
\newcommand{\bdry}{{\partial\Omega}}
\newcommand{\panel}{\Gamma}

\newcommand{\sglpot}{\mathcal{S}}
\newcommand{\dblpot}{\mathcal{D}}
\newcommand{\dVy}{\dif V_{\yy}}
\newcommand{\dSy}{\dif S_{\yy}}

\newcommand{\deltat}{{\delta t}}
\newcommand{\presvec}{\phi}

\newcommand{\conj}[1]{\bar{#1}}

\newcommand{\Npanels}{{N_p}}
\newcommand{\Nbdry}{{N_\bdry}}
\newcommand{\Ndom}{{N_\domain}}

\newcommand{\eqr}{\eqref}

\begin{document}
\hypersetup{pdfborder = 0 0 0, colorlinks=false,
 linkcolor=black,citecolor=black, filecolor=black, urlcolor=black, }

\begin{frontmatter}

\title{A Fast Integral Equation Method for the Two-Dimensional Navier-Stokes Equations}

\author[SFU, cor]{Ludvig af Klinteberg}
\author[NJIT]{Travis Askham}
\author[SFU]{Mary Catherine Kropinski}

\cortext[cor]{Corresponding author. E-mail address: \href{mailto:jafklint@sfu.ca}{jafklint@sfu.ca} }

\address[SFU]{Department of Mathematics, Simon Fraser University,
Burnaby, BC, Canada.}
\address[NJIT]{Department of Mathematical Sciences, New
  Jersey Institute of Technology, Newark, NJ, USA.}

\begin{abstract}
  % [MOTIVATION]  
  The integral equation approach to partial differential
  equations (PDEs) provides significant advantages in the
  numerical solution of the incompressible Navier-Stokes
  equations.
  In particular, the divergence-free condition and boundary
  conditions are handled naturally, and the ill-conditioning
  caused by high order terms in the PDE is preconditioned
  analytically.
  Despite these advantages, the adoption of integral equation
  methods has been slow due to a number of difficulties
  in their implementation.
  % [METHOD]
  This work describes a complete integral equation-based flow solver
  that builds on recently developed methods for singular quadrature
  and the solution of PDEs on complex domains, in combination with
  several more well-established numerical methods.
  % [KEY RESULTS]
  We apply this solver to flow problems on a number of geometries,
  both simple and challenging, studying its convergence properties and
  computational performance.
  % [CONCLUSIONS]
  This serves as a demonstration that it is now relatively
  straightforward to develop a robust, efficient, and flexible
  Navier-Stokes solver, using integral equation methods.
\end{abstract}


\end{frontmatter}

%%%%%%%%%%%%%%%%%%%%%%%%%%%%%%%%%%%%%%%%%%%%%%%%%%%%%%%%%%%%%%%%%%%%%%%%
\section{Introduction}
Fast integral equation methods (FIEMs) - techniques based on integral equations coupled with associated fast algorithms for the integral operators such as the Fast Multipole Method (FMM) \cite{Greengard1997} or fast direct solvers \cite{Martinsson2005} - have become the tool of choice for solving the linear, elliptic boundary value problems associated with the Laplace and Helmholtz equations, the equations of elasticity, the Stokes equation, and other classical elliptic equations in mathematical physics. For these boundary value problems, FIEMs have been a game changer: superior stability and accuracy, efficiency through dimension reduction and acceleration, and ease of adaptivity are a few of the more significant advantages.   

While nonlinear equations are not directly amenable to solution via integral equation methods, for some equations a suitable temporal discretization and/or linearization gives rise to a sequence of ``building-block'' partial differential equations (PDEs) of the form of the classic equations mentioned above. These can then be targeted with suitably-chosen FIEMs. The focus of this paper is to discuss a FIEM-based solver for the two-dimensional incompressible Navier-Stokes equations (INSE),
\begin{align}
  \begin{aligned}
    \pderiv{\uu}{t} + (\uu \cdot \nabla) \uu &= 
    -\nabla p + \frac{1}{\Re}\Delta \uu, && \xx \in \domain , \\
    \nabla \cdot \uu &= 0, && \xx \in \domain , \\
    \uu &= \ff, && \xx \in \bdry ,
  \end{aligned}
                   \label{eq:navier-stokes}
\end{align}
on a general bounded domain $\Omega$ and satisfying the no-slip boundary conditions. 
In this context, the basic idea is as follows. Applying a semi-implicit temporal discretization to the momentum equation in~\eqr{eq:navier-stokes} gives rise to a sequence of modified biharmonic (in the stream function formulation, c.f.~\cite{Jiang2013}) or modified Stokes equations, 
\begin{equation}
  \label{eqn:modstokes}
  \begin{aligned}
    (\alpha^{2} - \Delta)\uu + \Re \nabla p &= \FF,
      \quad &&\xx \in \domain, \\
    \nabla \cdot \uu &= 0, &&\xx \in \domain, \\
    \uu &= \ff, && \xx \in \bdry, 
  \end{aligned}
\end{equation}
where $\alpha^2 = O( \Re / \delta t)$, $\delta t$ is the time step, and with $\FF$ being updated at each time step. This is the so-called Rothe's method or Method of Lines Transpose. It is \eqr{eqn:modstokes}, then, that are recast by representing the solution as the sum of suitably chosen layer and volume potentials. This representation will satisfy the incompressibility constraint by construction, thereby eliminating the need for artificial boundary conditions associated with projection and other methods based on direct discretization of the Navier-Stokes equations. The focus of this paper, then, is to present a FIEM for \eqref{eqn:modstokes} in general, two-dimensional bounded domains and to examine their performance for a range of $\Re$.

This general approach was examined in \cite{Greengard1998} for the special case of flow inside a circular cylinder. It was demonstrated in this paper that FIEMs have a great deal of potential computing solutions for Reynolds number up to $O(1000)$ with both high spatial and temporal accuracy. Despite inspiring some other integral equation-based solvers for the INSE \cite{biros2002embedded}, this work has yet to be fully realized in general 2-D or 3-D domains. As outlined in \cite{Kropinski2011a}, there have been two critical stumbling blocks in achieving this goal. The first is how to evaluate the layer potentials close to the boundary where the kernels become nearly singular. There are many high-order quadrature methods that are well-suited for discretizing the integral operators and evaluating potentials away from the boundary \cite{Hao2014},  but these break down when the target becomes close to source points on the boundary. The second stumbling block is how to solve inhomogenous equations on general domains. A potential theoretic approach expresses the particular solution as a volume potential: a convolution of the fundamental solution with the right hand side of the equation. Unfortunately, available fast methods for evaluating layer potentials are generally in the form of ``box codes" such as \cite{Cheng2006} or the FFT. These require knowing the source data throughout a square or cube. In order for these box codes to be effective, this data must be of sufficient regularity throughout the box. For a boundary value problem on a general domain, this then requires extending the right-hand side with the required regularity throughout the bounding box.   

Recently, significant progress has been made on both the quadrature and function extension fronts.  Developing specialized quadrature that can be used to evaluate potentials arbitrarily close to the boundary and which employs the governing hierarchical acceleration strategy has been an area of active research.  One approach is ``Quadrature By Expansion'', or QBX \cite{Barnett2014, Klockner2013}, which uses local expansions with centres close to the boundary. QBX has the advantage that it can be extended to three dimensions but it can be challenging to integrate it effectively within an FMM. Another approach is an explicit kernel-split, panel-based Nystr\"{o}m method outlined in~\cite{Helsing2015, Helsing2008}. Here, the kernel is split into its smooth and regular parts, with the smooth parts being evaluated using a composite n-point Gaussian quadrature scheme, and the singular parts using specialized quadrature based on product integration.  In our case, the kernel-splitting of the Stokes kernel is considerably more complex with an additional complication that the underlying kernels require more spatial resolution for larger values of $\alpha$ (corresponding to large $\Re$ or small time steps). A straightforward implementation of \cite{Helsing2015} shows a breakdown in accuracy which requires careful local refinement to alleviate this \cite{AfKlinteberg2019}. The result is a quadrature scheme which is robust for arbitrary target points and values of $\alpha$. 

Two approaches have been investigated recently which extend functions throughout a bounding box in which the geometry of the problem is embedded. One approach uses extension via layer potentials \cite{Askham2016,Askham2017}.  The layer potential is determined by enforcing continuity across the boundary of the domain through solving an integral equation, and it is then used to construct the function extension in the complement of the domain. Requiring a higher degree of continuity means the integral equation is based on solving a higher-order partial differential equation and this can get prohibitively complex. A second approach, which is the approach we take here, is to smoothly extend the right-hand-side of \eqref{eqn:modstokes_particular}  via extrapolation using radial basis functions~\cite{Fryklund2018,Fryklund2019}. This method is called Partition of Unity eXtension, or PUX. 
It is then simple to convolve this extended function with the periodic Green's function for these equations using a fast Fourier transform (FFT)
and the solution may be evaluated on the
boundary rapidly with a non-uniform FFT
\cite{Barnett2018}.

The methods we propose in this paper consist of coupling together existing and custom-developed state-of-the-art tools for constructing and evaluating the layer and volume potentials. The key elements of the overall approach and the road map to where these elements are discussed in this paper are summarized below:
\begin{enumerate}
\item Using PUX, the function $\FF$ is extended onto a uniform grid discretizing the bounding box. The volume potential is evaluated on this grid via an FFT and evaluated on $\bdry$ via a non-uniform FFT. This is discussed in section 3.1.
\item To the volume potential, layer potentials are added which ensure the solution satisfies the no-slip boundary conditions $\uu = \ff$. This involves solving a system of Fredholm integral equations of the second kind for the homogenous modified Stokes equations (see section 2.2). After applying quadrature (see section 4), the resulting linear system is dense and of the form of $[I + K] \mmu = {\bf b}^N$, where $N$ represents the time step and $K$ is has low-rank structures. Since the boundary $\bdry$ is fixed, the system matrix $[I+K]$ does not change throughout time. We use a fast direct method \cite{Gillman2012, Marple2016} which exploits the low-rank off-diagonal blocks to invert the system matrix as a precomputation. This inverse can be applied to the updated vector ${\bf b}^N$ at each time step at a fraction of the time. All computations are linear in computational time. This is discussed in section 3.2.2. 
\item Given the density of the layer potentials, $\mmu$, from the previous step,  $\uu$ and $\nabla \uu$ are evaluated at all regular grid points in $\Omega$. Special-purpose quadrature is used when evaluating the layer potentials close to the source terms \cite{AfKlinteberg2019} (see section 4). These point-to-point interactions are accelerated by an FMM \cite{Greengard1997} based on the separation of variables formulae derived in \cite{Askham2018}, which are stable for all values of $\alpha$ (section 3.2.1). These values are used to construct a new right hand side $\FF$ in $\Omega$. 
\item The above steps are repeated for the next time step. The solution is advanced using a semi-implicit spectral deferred correction method \cite{Minion2003}. This is discussed in section 3.3. 
\end{enumerate}
The result is a Navier-Stokes solver which has linear or near-linear scaling in computational time at each time step, specifically the computational cost is $\mathcal O( \Ndom\log\Ndom + \Ndom + \Nbdry)$, where $\Ndom$ and
$\Nbdry$ denote the number of discretization nodes in the domain and
on the boundary, respectively. We show in section 5.5 that the stability requirement is $\delta t = O(1/\Re)$. 
The numerical experiments, which we present in section 5, indicate that our methods are well suited for problems with moderate values of the Reynolds number. 

Our approach could be described as an embedded boundary approach, which is a popular approach for methods based on direct discretization of the PDE, such as the Immersed Boundary (IB) method of
Peskin \cite{Peskin2002a}. Here,  the effect of the boundary is smeared to
neighboring grid cells (including the cut cells)
by using approximations to the Dirac $\delta$
function.
The original IB framework was limited to
low order accuracy; recently, Stein et al.
\cite{Stein2017} presented
a higher-order IB approach based on envisioning
a smooth extension of the solution outside
of the original domain.
The ideas behind the resulting method, IB Smooth
Extension (IBSE), bear some resemblance to the
present work but are based on introducing new
unknowns for the solution outside of the domain
as opposed to extending
the known inhomogeneity, as we describe below.
Further references for embedded boundary
methods in the PDE literature can be
found in \cite{Stein2017}.

There has been recent work on time-dependent integral
equation methods for the unsteady Stokes equations
\cite{Greengard2018,Greengard2018a,Wang2018}, which could form the
basis of a FIEM-based Navier-Stokes solver quite different from the
one presented here. 
These methods are still in the early stages of development, but are
nevertheless promising. In particular, their formulation allows for a
natural treatment of moving boundaries, and would likely have
advantages in terms of stability.

The remainder of this paper is organized as follows. 
\Cref{sec:formulation} introduces the underlying mathematical formulation.
\Cref{sec:numerical-methods} outlines the framework of numerical methods.
\Cref{sec:quadrature} gives a detailed description of the quadrature used for singular and nearly singular integrals.
\Cref{sec:numerical-examples} reports on results from numerical experiments.
Finally, concluding remarks are provided in \cref{sec:conclusion}, while supplementary material is available in 
\cref{sec:layer-potentials,sec:fundamental-solutions,sec:power-series}.

%%% Local Variables: 
%%% TeX-master: "nseIE"
%%% End: 

%%%%%%%%%%%%%%%%%%%%%%%%%%%%%%%%%%%%%%%%%%%%%%%%%%%%%%%%%%%%%%%%%%%%%%%%
\section{Formulation}
\label{sec:formulation}

Many numerical approaches to solving the Navier-Stokes
equations begin by performing a discretization of
the equations in time.
Because it is simpler to handle the nonlinear term
explicitly but less restrictive to handle the stiff
diffusion term implicitly, it is common to use an
implicit-explicit (IMEX) method \cite{Ascher1995}.
Using the first-order IMEX-Euler rule results in the
time-independent elliptic PDE
\begin{equation}
  \label{eq:navier-stokes-imex}
\begin{aligned}
  \frac{\uu^{N+1} - \uu^N}{\deltat} - \frac{1}{\Re}\Delta \uu^{N+1} +
  \nabla p^{N+1} &= -(\uu^N \cdot \nabla) \uu^N, && \xx \in \domain,  \\
  \nabla \cdot \uu^{N+1} &= 0, && \xx \in \domain, \\
  \uu^{N+1} &= \ff, && \xx \in \bdry.
\end{aligned}
\end{equation}
The higher-order IMEX rules result in qualitatively similar
elliptic PDEs.
Further, the spectral deferred corrections (SDC)
method \cite{Dutt2000} can be used to combine many
low order IMEX steps into a high order scheme
for advection-diffusion problems \cite{Minion2003},
which will be the approach of this paper.

It is in the solution of the elliptic PDE
\eqref{eq:navier-stokes-imex}, called the modified Stokes
equations\footnote{These equations may also be referred to as the
  linearized unsteady Stokes equations \cite{Pozrikidis1989}, or the
  Brinkman equations \cite{Brinkman1947}.}
where the PDE discretization and integral equation methods diverge.
Because the equation is linear in the unknowns
$\uu^{N+1}$, the solution of~\eqref{eq:navier-stokes-imex} can
be split into a particular and a homogeneous problem,
\begin{equation}
  \label{eqn:modstokes_particular}
  \begin{aligned}
    (\alpha^{2} - \Delta)\uu^P + \Re \nabla p^P &= \FF,
      \quad &&\xx \in \domain, \\
    \nabla \cdot \uu^P &= 0, &&\xx \in \domain,
  \end{aligned}
\end{equation}
and
\begin{equation}
  \label{eqn:modstokes_homogeneous}
  \begin{aligned}
    (\alpha^{2} - \Delta)\uu^H + \Re \nabla p^H &= 0, 
      \quad &&\xx \in \domain, \\
    \nabla \cdot \uu^H &= 0, &&\xx \in \domain, \\
    \uu^H &= \ff - \uu^P, &&\xx \in \bdry \; ,
  \end{aligned}
\end{equation}
where
\begin{align}
  \alpha^2 = \frac{\Re}{\deltat},
  \quad \text{ and } \quad
  \FF = \alpha^2 \uu^N - \Re (\uu^N\cdot \nabla)\uu^N \label{eq:alpha_and_rhs} .
\end{align}
The forced PDE~\eqref{eqn:modstokes_particular} has no
boundary conditions, so that a candidate particular solution
may be computed by convolving the right-hand-side
of~\eqref{eqn:modstokes_particular} with a free-space
or periodic Green's function
--- in particular, a Green's function which is independent
of the domain and analytically known.
Regardless of the particular solution chosen, the
homogeneous equation~\eqref{eqn:modstokes_homogeneous}
is specified so that the composite solution,
\begin{align}
  \uu^{N+1} &= \uu^P + \uu^H, \\
  p^{N+1} &= p^P + p^H,
\end{align}
satisfies the original problem~\eqref{eq:navier-stokes-imex}.
For this reason, the homogeneous solution is sometimes referred to as a
boundary correction.

In an integral equation method, the homogeneous solution
is generally represented by the convolution of an unknown density
defined on the boundary, $\bdry$, with an integral
kernel given in terms of the equivalent of a charge or
a dipole of the equations~\eqref{eqn:modstokes_homogeneous},
which a priori satisfies the divergence-free condition and
the homogeneous equations inside the domain.
For a well-chosen representation, the enforcement
of the boundary conditions gives a well-conditioned
equation for the density on the boundary.
Once this boundary integral equation is solved,
the homogeneous solution may be evaluated at any point
in the domain by convolution.

The central equation of our method is the inhomogeneous modified
Stokes equations with Dirichlet boundary conditions,
\begin{equation}
  \label{eq:modstokes}
  \begin{aligned}
    (\alpha^{2} - \Delta)\uu + \nabla p &= \FF,
      \quad &&\xx \in \domain, \\
    \nabla \cdot \uu &= 0, &&\xx \in \domain, \\
    \uu &= \ff, &&\xx \in \bdry.
  \end{aligned}
\end{equation}
This stationary problem is the result of a semi-implicit
discretization of the Navier-Stokes equations
\eqref{eq:navier-stokes-imex}.
Because of the divergence-free condition, these equations,
like the Navier-Stokes equations, imply a compatibility
condition on $\ff$, namely

\begin{equation}
  \label{eq:fcompat}
  \int_\bdry \ff \cdot \v\nhat \, \dSy = 0 \; .
\end{equation}

In analogy with the regular Stokes
equations, we call the free-space Green's function of these equations
the modified stokeslet tensor, or simply the stokeslet. It is defined as
\begin{align}
  \label{eq:stokeslet}
  S(\v x, \v y) = (\nabla\otimes\nabla -\Delta) G(\v x, \v y),
\end{align}
where 
\begin{align} 
  G(\xx,\yy) = -\frac{1}{2\pi\alpha^2}
  \left (\log \|\xx - \yy\| +
  K_0(\alpha\|\xx - \yy\|) \right) \label{eq:modbiharm}
\end{align}
is the free-space Green's function of the modified biharmonic equation,
i.e. $\Delta(\Delta - \alpha^2)G = \delta$ \cite{Jiang2013}. The
corresponding pressure vector $\v\presvec$ and stress tensor $T$
(which we will refer to as the stresslet) are defined as
\begin{align}
  \presvec_i(\xx, \yy) &= (\Delta-\alpha^2) \dpd{G(\xx,\yy)}{x_i}, \\
  T_{ijk}(\xx, \yy) &= -\delta_{ik} \presvec_j(\xx, \yy) 
                        + \dpd{S_{ij}(\xx, \yy)}{x_k}  
                        + \dpd{S_{kj}(\xx, \yy)}{x_i} \label{eq:stresslet}.
\end{align}
Note that if $\uu = S \mmu$ is the field induced by a stokeslet
with charge $\mmu$, then the stress tensor for that field, $\ssigma$,
is given by
\begin{equation}
\sigma_{ik} = T_{ijk} \mu_j .
\end{equation}
Closed-form expressions for $S$, $\v\presvec$, and $T$ are derived in
\cref{sec:fundamental-solutions}.

Our solution strategy for the equations \eqref{eq:modstokes} is based
on two fundamental principles, which are typical of integral equation
methods.
The first principle is the use of potential theory, in which
the solution is represented using
convolutions with Green's functions of the equations. By construction,
this guarantees that the divergence-free condition is satisfied, in
contrast with methods based on direct discretization of the
PDE. The second principle is the use of a composite solution, which is
formed by the solutions to the particular problem
\eqref{eqn:modstokes_particular} and the homogeneous problem
\eqref{eqn:modstokes_homogeneous}. This split allows us to use
numerical methods that are fast and geometrically flexible, as we
shall see in \cref{sec:numerical-methods}. In the remainder
of this section, we describe in detail how the Green's functions
are used to solve the particular and homogeneous problems.

\subsection{Particular solution}
\label{sec:particular-solution}

From a mathematical point of view, the most straightforward solution to the
particular problem \eqref{eqn:modstokes_particular} is the volume
potential,
\begin{align}
  \uu^{P}(\xx) = \int_{\domain} S_{ij}(\xx,\yy) \FF_j(\yy) \dVy,
  \label{eq:vol_pot}
\end{align}
which solves \eqref{eqn:modstokes_particular} using free-space
boundary conditions. This is however not practical, as evaluating the
integral \eqref{eq:vol_pot} accurately over a complex domain would
present a significant numerical challenge. Instead, we shall for the
moment assume that $\domain$ is embedded in a square domain
$B=[0,2\pi]^{2}$ (see \cref{fig:embedded_domain}), and that there
exists a compactly supported function $\v F^e$ that is a q-continuous
extension of $\v F$ out into $B$,
\begin{align}
      &\FF^e = \FF \quad \text{in } \Omega,\\
      &\domain \subset \text{supp}\{\FF^{e}\} \subset B,\\
      &\FF^e \in C^q(B), \quad q > 0.
\end{align}
We will return to the computation of the extension $\v F^e$ in
\cref{sec:extension}. Having $\v F^e$ lets us compute $\uu^P$ as the
solution to the extended particular problem,
\begin{equation}
  \label{eqn:modstokes_particular_box}
  \begin{aligned}
    (\alpha^{2} - \Delta)\uu^P + \Re \nabla p^P &= \FF^e,
      \quad &&\xx \in B, \\
    \nabla \cdot \uu^P &= 0, &&\xx \in B .
  \end{aligned}
\end{equation}
Note that the solution to the extended particular problem
\eqref{eqn:modstokes_particular_box} also satisfies the original
particular problem \eqref{eqn:modstokes_particular}. Now, since $B$ is
a simple domain (a box), it would be straightforward to compute the
particular solution as the volume potential over $B$ with the extended
function $\v F^e$. Alternatively, and this is the path we will follow,
we can compute it using a Fourier series. This amounts to solving
\eqref{eqn:modstokes_particular_box} using periodic boundary
conditions, which we can do since $\FF^e$ is compactly supported, and
therefore also periodic. Let
\begin{align}
  u^P_j(\xx) &= \sum_{\kk \in \ZZ^{2}} \hat{u}_j(\kk)
               e^{i \kk \cdot \xx},
  &
    p^P(\xx) &= \sum_{\kk \in \ZZ^{2}} \hat p (\kk)
                           e^{i \kk \cdot \xx},
  &
    F^e_j(\xx) &= \sum_{\kk \in \ZZ^{2}} \hat F^e_j (\kk) 
                 e^{i \kk \cdot \xx}.
\end{align}
Replacement into \eqref{eqn:modstokes_particular_box} and matching of
coefficients gives the linear system of equations
\begin{align}
  (\alpha^2 + |\kk|^2) \hat u_j(\kk) + \Re i k_j \hat p(\kk) &= \hat F^e_j(\kk), \\
  i k_j \hat u_j(\kk) & = 0,
\end{align}
with solution
% \begin{align}
%   % Derivation:
%   &- \Re |\kk|^2 \hat p(\kk) = i k_j \hat F^e_j(\kk) \\
%   &\Re \hat p(\kk) = -\frac{i k_l \hat F^e_l(\kk)}{|\kk|^2} \\
%   &(\alpha^2 + |\kk|^2) \hat u_j(\kk) - i k_j \frac{i k_l \hat F^e_l(\kk)}{|\kk|^2} = \hat F^e_j(\kk),  \\
%   &\hat u_j(\kk) = \frac{\hat F^e_j(\kk) - k_j \frac{k_l \hat F^e_l(\kk)}{|\kk|^2}}{\alpha^2 + |\kk|^2},  
% \end{align}
\begin{align}
  \begin{aligned}
    \hat u_j(\kk)
    &= \frac{ |\kk|^2\delta_{jl} - k_jk_l }{|\kk|^2(\alpha^2 + |\kk|^2)} \hat F^e_l(\kk), \quad \v k \ne 0,\\
    \hat u_j(0) &= \frac{1}{\alpha^2} \hat F^e_j(0).    
  \end{aligned}
  \label{eq:uhat}
\end{align}
To rephrase the above, given a smooth extension of $\FF$
with Fourier coefficients $\hat\FF^e$, a solution
of the particular equation, \cref{eqn:modstokes_particular},
is given by the convolution
\begin{align}
  u^P_j(\xx) &= \sum_{\kk \in \ZZ^{2}} \hat S_{jl}(\kk) \hat F^e_l(\kk)
               e^{i \kk \cdot \xx} \label{eq:uP_Fourier}, \quad \xx \in \domain,
\end{align}
where $\hat S$ is the Fourier multiplier corresponding
to the periodic stokeslet, i.e.
\begin{align}
  \begin{aligned}
    \hat S_{jl}(\kk) &= \frac{ \delta_{jl} - \hat k_j \hat k_l}{\alpha^2 + |\kk|^2},
    \quad \v k \ne 0, \quad \v{\hat k} = \kk / |\kk|, \\
    \hat S_{jl}(0) &= \frac{\delta_{jl}}{\alpha^2} .\\
  \end{aligned}
  \label{eq:periodic_stokeslet}
\end{align}
We also need the gradient of the solution, in order to compute the
right hand side of the next time step \eqref{eq:alpha_and_rhs}. This is straightforward to
compute using Fourier differentiation,
\begin{align}
  \dpd{u^P_j(\xx)}{x_m} &= \sum_{\kk \in \ZZ^{2}} i k_m \hat S_{jl}(\kk) \hat F^e_l(\kk)
               e^{i \kk \cdot \xx} \label{eq:uP_Fourier_deriv}, \quad \xx \in \domain.
\end{align}
This solution procedure has some nice properties.
First, once the extension $\FF^e$ is available,
computing the particular solution does not require
that a linear system be solved; instead, the computation
is based on evaluating a formula for the solution.
Further, for smaller values of $\alpha$, the multiplier
$\hat S$ is strongly smoothing, i.e. it decreases the
magnitude of higher Fourier modes.
This property is both consistent with the physics and
stabilizing from a numerical point of view.
Unfortunately, this smoothing property is reduced for large values of
$\alpha$, corresponding to high Reynolds number flow.

\subsection{Homogeneous solution}
\label{sec:homogeneous-solution}

The solution to the homogeneous problem
\eqref{eqn:modstokes_homogeneous} can be represented by a double layer
potential, which we denote by
\begin{align}
  \uu^{H}(\xx) = \dblpot[\mmu](\xx), \quad \xx \in \domain. \label{eq:uH_dbl}
\end{align}
For a given double layer density $\mmu$ defined on $\bdry$, the
potential is defined as the boundary integral
\begin{align}
  \dblpot_i[\mmu](\xx) &= \int_{\bdry} D_{ij}(\xx,\yy) \mu_j(\yy) \dSy, \label{eq:dbl_lyr_pot}
\end{align}
where the double layer kernel, $D$, is defined using the stresslet
\eqref{eq:stresslet} and the surface unit normal $\v \nhat(\v y)$
pointing into $\domain$,
\begin{align}
  D_{ij}(\v x, \v y) = -T_{jik}(\xx, \yy) \nhat_k(\v y) .
  \label{eq:dbl_lyr_ker}
\end{align}
By construction, this layer potential is divergence-free
and satisfies the homogeneous modified Stokes equation
in $\domain$.

We will now briefly describe how a double layer potential
is used to solve the homogeneous Dirichlet problem.
Much of the following depends on a number of
technical results concerning the double layer potential
which can be found in the literature
\cite{KimSangtae1991M:pa,Pozrikidis1992,biros2002embedded} and
are included in \Cref{sec:layer-potentials} for
convenience.
In contrast with the Stokes equation, we note that all of
the following holds without modification for multiply-connected
$\domain$.
The appropriate density $\mmu$ is determined by enforcing
the Dirichlet boundary condition $\uu^{H}=\ff-\uu^P$ in the limit
$\Omega \ni \xx \to \bdry$.
Using the jump conditions for modified Stokes layer potentials (see
\Cref{lemma:jump-conds} in \cref{sec:layer-potentials}), we obtain a
second-kind integral equation in $\mmu$,
\begin{align}
  \frac{1}{2}\mmu + \DD[\mmu] = 
  \ff - \uu^P, \quad \xx \in \bdry .
  \label{eqn:BIE_nocorr}
\end{align}
This integral equation has a null space of rank 1
(see \Cref{lemma:nullspaces} in \cref{sec:layer-potentials}), so in order to have a unique
solution we add the term
\begin{align}
  \mathcal W[\mmu](\xx)
  = \frac{\v\nhat(\xx)}{\int_{\bdry} \dif S} \int_{\bdry}
  \mmu(\yy) \cdot \v\nhat(\yy) \dSy \; ,
\end{align}
which is sometimes referred to as a nullspace correction.
Then, the homogeneous solution is given by a double layer potential
\eqref{eq:uH_dbl}, with the density satisfying
\begin{align}
  \frac{1}{2}\mmu + \DD[\mmu] + \mathcal W[\mmu] = 
  \ff - \uu^P, \quad \xx \in \bdry .
  \label{eqn:BIE}
\end{align}
The new equation, \cref{eqn:BIE}, is invertible for any
right hand side and the term $\mathcal{W} [\mmu]$ is
zero provided that the right hand side satisfies the
condition (see \Cref{lemma:invertibility} in \cref{sec:layer-potentials})

\begin{equation}
  \int_\bdry \left ( \ff - \uu^p \right ) \cdot \v\nhat
  \, \dSy = 0 \; .
\end{equation}
This automatically holds for $\ff$ satisfying
the compatibility condition, \cref{eq:fcompat},
and any particular solution, $\uu^p$.
Once this equation is solved, $\uu^H$ can be evaluated anywhere inside
$\Omega$ using \eqref{eq:dbl_lyr_pot}.

%%% Local Variables: 
%%% TeX-master: "nseIE"
%%% End: 

%%%%%%%%%%%%%%%%%%%%%%%%%%%%%%%%%%%%%%%%%%%%%%%%%%%%%%%%%%%%%%%%%%%%%%%%
\section{Numerical methods}
\label{sec:numerical-methods}

This section gives a complete outline of the numerical methods
used. Quadrature will only be briefly discussed in this section, and
more thoroughly covered in \cref{sec:quadrature}.

\begin{figure}
  \centering
  \includegraphics[width=0.25\textwidth]{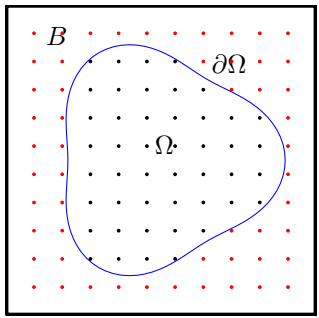}
  \caption{An illustration of the embedded boundary approach. The
    boundary of the domain is ``floating'' on top of a uniform grid,
    defined on a bounding box $B$. The forcing function $\FF$,
    initially only known at grid points inside $\domain$ (marked
    black), is smoothly extended to the function $\FF^e$, which is
    defined at all grid points in $B$.}
  \label{fig:embedded_domain}
\end{figure}

%%%%%%%%%%%%%%%%%%%%%%%%%%%%%%%%%%%%%%%%%%%%%%%%%%%%%%%%%%%%%%%%%%%%%%%%
\subsection{Evaluating the particular solution}

%%%%%%%%%%%%%%%%%%%%%%%%%%%%%%%%%%%%%%%%%%%%%%%%%%%%%%%%%%%%%%%%%%%%%%%%
\subsubsection{Smooth function extensions}
\label{sec:extension}

In order compute the particular solution $\uu^P$ using the Fourier
method of \cref{sec:particular-solution}, we need an extended forcing
function $\FF^e$ with a high order of regularity. To compute such a
function extension, we use the method known as PUX, for Partition of
Unity eXtension, introduced in \cite{Fryklund2018}. The
method is based on overlapping circular partitions with their centers
distributed along the boundary $\bdry$. Inside each partition radial
basis functions (RBFs) are defined, which interpolate $\FF$ at each
point of the partition inside $\domain$, giving a local extrapolation
of $\FF$ to the other points in the partition. These exterior values
in adjacent partitions, which can be multiply-defined due to the
overlap, are then blended using a high-order partition of unity (POU)
function, forming a function extension with high regularity. Finally,
the extension is given compact support by blending it with a second
layer of ``zero partitions'', located just outside $\domain$.

The current version of PUX is limited to a uniform grid in $B$, and we
will use a uniform grid for that reason. Given a grid and a boundary,
\cite{Fryklund2018} provides heuristics for determining the parameters
of the method; the only free parameters are the RBF shape parameter
$\epsilon$ (which is set to 2 throughout) and the partition radius
$R$. The accuracy largely depends on the number of grid points spanned
by $R$; the more grid points spanned, the better the decay of the
extension can be resolved, up to a limit. The range of points spanned
by $R$ is typically between 20 and 50.

%%%%%%%%%%%%%%%%%%%%%%%%%%%%%%%%%%%%%%%%%%%%%%%%%%%%%%%%%%%%%%%%%%%%%%%%
\subsubsection{Evaluating the volume potential}
\label{sec:volumeFourier}

After applying the PUX algorithm, we have the values of the smooth
extension $\FF^e$ on a uniform grid in $B$. It is then straightforward
to compute the particular solution $\uu^P$ \eqref{eq:uP_Fourier} in a
fast way, using fast Fourier transforms (FFTs). First, two FFTs give
us the Fourier coefficients $\hat\FF^e(\kk)$. We can then compute the
solution using two inverse FFTs for the volume grid, and two
non-uniform FFTs of type 2 \cite{Barnett2018} for the discretization
points on $\bdry$. The gradient \eqref{eq:uP_Fourier_deriv} is then
computed using four additional inverse FFTs. Since the gradient is
only needed on the volume grid, it does not require any non-uniform
FFTs.

For a discretization with $N_B$ uniform points in B, and $N_\bdry$
nonuniform points on $\bdry$, the computational cost is
$\bigO (N_B \log N_B + N_\bdry)$.  The method has spectral accuracy if
$\FF$ and its extension are smooth. If they are not smooth, the order
of accuracy is related to the regularity of the extended function.

%%%%%%%%%%%%%%%%%%%%%%%%%%%%%%%%%%%%%%%%%%%%%%%%%%%%%%%%%%%%%%%%%
\subsection{Computing the homogeneous correction}
\label{sec:comp-homog-corr}

In order to solve the integral equation \eqref{eqn:BIE}, we first need
to discretize the boundary $\bdry$ using a suitable quadrature. For
this, we use a composite Gauss-Legendre quadrature, wherein $\bdry$ is
subdivided into $\Npanels$ panels $\Gamma_i$ of equal arclength,
\begin{align}
  \bdry = \bigcup_{i=1}^\Npanels \Gamma_i,
\end{align}
each of which is discretized using an $n$-point Gauss-Legendre
quadrature (we use $n=16$ throughout this work). If we let $\v y_i(t)$ be
the mapping from the standard interval $[-1,1]$ to the panel
$\Gamma_i$, then
\begin{align}
  \int_{\Gamma} f\left(\v y\right) \dSy
  = \sum_{i=1}^\Npanels \int_{\Gamma_i} f\left(\v y\right) \dSy 
  = \sum_{i=1}^\Npanels \int_{-1}^1 f(\v y_i(t)) \abs{\v y_i'(t)} \dif t
  \approx \sum_{i=1}^\Npanels \sum_{j=1}^n f(\v y_i(t_j)) \abs{\v y_i'(t_j)} \lambda_j .
\end{align}
Here $(t_j,\lambda_j)_{j=1}^n$ are the standard Gauss-Legendre nodes
and weights on $[-1,1]$. For convenience, we now combine the double
summation over $i$ and $j$ into a single index $m$ that runs to
$N_\bdry = n \Npanels$, and let $\xx_m = \v y_i(t_j)$ and
$w_m = \abs{\v y_i'(t_j)} \lambda_j$. This gives us a set of nodes and
weights $(\xx_m, w_m)_{m=1}^{N_\bdry}$ for approximating integrals
over $\bdry$. Given the double layer density $\mmu$ at the nodes, we
can then approximate the double layer potential \eqref{eq:dbl_lyr_pot}
anywhere in $\domain$ as
\begin{align}
  \tilde\dblpot_i[\mmu](\xx) &= \sum_{m=1}^{\Nbdry} D_{ij}(\xx,\xx_m) \mu_j(\xx_m) w_m,
  \quad \xx \in \domain. \label{eq:dbl_lyr_quad}
\end{align}
As long as $\bdry$ and $\mmu$ are properly resolved, this quadrature
will be accurate for target points $\xx$ that are well-separated from
$\bdry$. However, straight Gauss-Legendre quadrature will not be
accurate for target points close to, or on, the boundary $\bdry$. This
is due to the singularities in the stresslet \eqref{eq:stresslet}. We
will discuss specialized quadrature methods for such cases in
\cref{sec:quadrature}. For now, we write $w_m = w_m(\xx)$, with the
understanding that these weights are only dependent on $\xx$ at the
finite number of boundary points $\xx_m$ located close to $\xx$.

We approximate the integral equation using the quadrature
\eqref{eq:dbl_lyr_quad}, including specialized quadrature weights
where needed, and enforce it at the quadrature nodes (this is known as
the Nystr\"om method). This gives us the $2\Nbdry \times 2\Nbdry$
dense linear system
\begin{align}
  \frac{1}{2}\mmu(\xx_n) + \sum_{m=1}^{\Nbdry}
  \left(
  D(\xx_n,\xx_m)
  + \v\nhat(\xx_n) \v\nhat(\xx_m)^T
  \right) \mmu(\xx_m) w_m(\v x_n)
  = \ff(\xx_n) - \uu^P(\xx_n), \quad n=1,\ldots,\Nbdry.
  \label{eq:BIE_quad}
\end{align}
As this is a second-kind integral equation, it
can be solved iteratively using GMRES \cite{Saad1986} in a number of
a iterations that is bounded as $\Nbdry$ grows. We can then evaluate
the homogeneous solution on the domain grid using
\eqref{eq:dbl_lyr_quad}. If done naively, both of these steps have
quadratic costs, which is infeasible for large systems. To get around
this, we make use of fast hierarchical methods, described in
\cref{sec:fast-summation,sec:fast-direct-solver}.

%%%%%%%%%%%%%%%%%%%%%%%%%%%%%%%%%%%%%%%%%%%%%%%%%%%%%%%%%%%%%%%%%%%%%%%%%%%%%
\subsubsection{Fast summation}
\label{sec:fast-summation}

Let $\Ndom\leq N_B$ denote the number of grid points which lie inside
$\domain$ (in \Cref{fig:embedded_domain} this is the number of black dots).
A naive evaluation of \cref{eq:dbl_lyr_quad} at each
of these points would require $O(\Ndom\Nbdry)$ operations. Using the
fast multipole method (FMM) \cite{Greengard1997}, it is possible to
reduce this computational cost to $O(\Ndom + \Nbdry)$, by taking
advantage of special structure in the sum. In this section, we describe
some components of an FMM based on the special functions defined
in \cite{Askham2018}.

\paragraph{Nomenclature}

In analogy with physics, we refer to a Green's function $G(\xx,\yy)$
as a charge located at $\yy$, a derivative of
$G$ as a dipole, a second derivative of $G$ as a quadrupole,
and a third derivative of $G$ as an octopole. By the definition of
the double layer potential, the quadrature formula \cref{eq:dbl_lyr_quad}
is a sum of dipoles and octopoles for the modified biharmonic
equation which are located along $\bdry$. These
charges and higher moments are referred to as sources and the
points at which we would like to evaluate the field are called
targets.

\paragraph{Well-separated interactions}

\begin{figure*}[h]
  \centering
  \subfloat[][]{
    \includegraphics[height=0.3\textwidth]{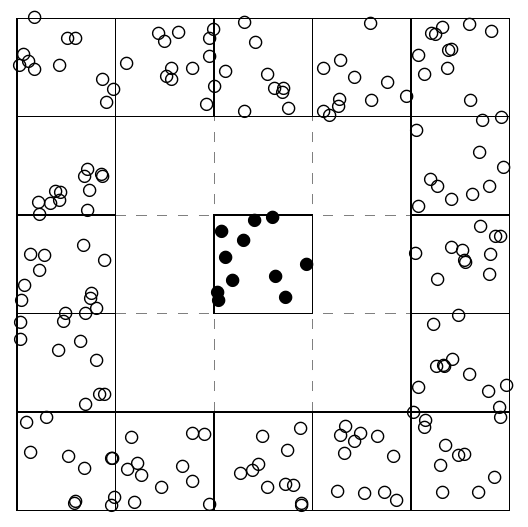}
    \label{fig:well-sep}
  }
  \subfloat[][]{
    \includegraphics[height=0.3\textwidth]{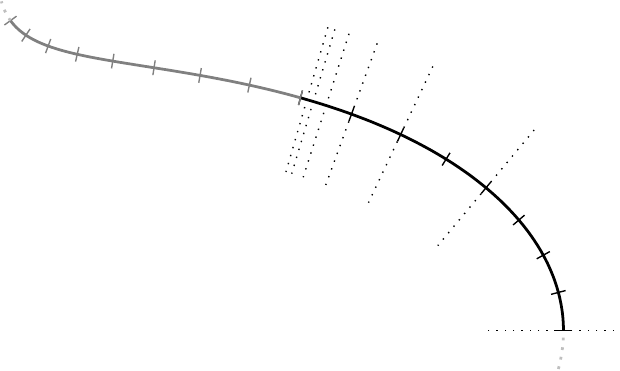}
    \label{fig:lowrank_adjacent_1d}
  }
  
  \caption{(a) We plot a collection of white dots
    which are well-separated from the black dots.
    (b) An illustration of the property that adjacent pieces of
    a one dimensional curve are of relatively low rank. The
    points on the gray section of the curve are well-separated
    from each set of points between the dotted separators on
    the black section of the curve.}
\end{figure*}

Suppose you can define two boxes around a set of points
such that they share the same center, the larger box
has three times the side length, and the smaller box
contains all of the points. Points located outside of the
larger box are said to be well-separated from the original
points, see \cref{fig:well-sep} for an illustration.
An important property of well-separated interactions is
that the field induced by a set of source points is
smooth from the perspective of well-separated points
and can be approximated accurately using relatively few
terms in a separation of variables expansion, 
independent of the number of source points.

\paragraph{Separation-of-variables representations}

Because the field induced by the sources
satisfies the modified biharmonic equation,
standard separation-of-variables
techniques apply. Let $I_n(\alpha r)$ and $K_n(\alpha r)$
denote the modified Bessel functions of the first and second
kind, respectively. Classically, a field which satisfies
the modified biharmonic equation can be represented in
the exterior of the disc as a series in $\log r$, $r^{-|n|}$,
and $K_n(\alpha r)$ and in the interior of a disc as a
series in $r^{|n|}$ and $I_n(\alpha r)$, where $r$ is the
radial part of polar coordinates centered in the disc.
Using these classical functions directly leads to catastrophic
cancellation when $\alpha r$ is small. 
Following \cite{Askham2018}, we define
the functions $Q_n$ and $P_n$ via

\begin{align}
  \label{eq:Qn}
  Q_n(r) &= K_n(\alpha r) - \dfrac{2^{|n|-1}
    \left (|n|-1\right )!}{\alpha^{|n|}r^{|n|}} \; , \; |n| > 0 \; , \\
  Q_0(r) &= K_0(\alpha r) + \log(r) \; ,  \label{eq:Q0}
\end{align}
and

\begin{equation} \label{eq:Pn}
  P_n(r) = I_n(\alpha r) 
- \left ( \dfrac{\alpha r}{2} \right )^{|n|} \dfrac1{|n|!} \; .
\end{equation}

Consider a disc centered at a point $\yy$ and
let $(r,\theta)$ denote the polar coordinates
of a point $\xx-\yy$.
Then, the field, $\psi$, induced by a collection of
sources inside the disc can be represented by 

\begin{equation}
  \psi(\xx) = \sum_{n=-\infty}^\infty \left (\alpha_n Q_n( r)
  +\beta_n K_n(\alpha r) \right) e^{in\theta} \; , \label{eq:mpole}
\end{equation}
which we will refer to as a multipole expansion,
and the field, $\phi$, induced by a collection of sources
outside the disc can be represented by

\begin{equation}
  \phi(\xx) = \sum_{n=-\infty}^\infty \left (a_n P_n(\alpha r) +
  b_n r^{|n|} \right ) e^{in\theta} \; , \label{eq:loc}
\end{equation}
which we will refer to as a local expansion.
These expansions can be manipulated stably.
In particular, \cite{Askham2018} provides
stable formulas for shifting the centers of
the expansions and converting a multipole
expansion about one center into a local expansion
about another. When evaluating the expansions
for a collection of sources which are well separated
from the targets, the standard decay rates in
\cite{greengard-1987,Cheng2006} apply, so that
truncating the sums in \cref{eq:mpole,eq:loc}
to order $\pfmm=40$, i.e. retaining only the terms
with $|n| \leq \pfmm$, is sufficient to
achieve 12 digits of accuracy.

\paragraph{Tree structure}

To deal with arbitrary collections
of sources and targets, as opposed to well-separated ones,
the source and target points are arranged in a
hierarchical quad-tree structure in the FMM. Starting from
a root box containing all sources and targets, the tree
nodes (boxes) are subdivided equally into four quadrants,
which are its children nodes, until they have at most a
prescribed total number of source and target points inside.
Typically, the maximum number
of points is set to the order of the expansions used
to represent the well-separated interactions.
Many FMM algorithms are written for arbitrary quad-trees;
for simplicity, we add a level-restriction to
the tree: all adjacent leaf boxes should be within
one level of refinement.

The FMM gains its efficiency by organizing
the interactions between sources and targets
hierarchically in this tree so that most interactions
are treated as well-separated. This is done using
a number of standard techniques, see
\cite{greengard-1987,Greengard1997} for a
reference.
Below, we only describe how multipole and
local expansions are formed on leaf boxes,
as the other techniques are a straightforward
application of the translation formulas in
\cite{Askham2018}.

\paragraph{Form multipole and form local}

Consider the sources contained within a leaf
box. Note that a charge can be represented
perfectly as a 0-term expansion about its own
center because $G(\xx,\yy) = -Q_0(r)/(2\pi\alpha^2)$,
with $r = |\xx-\yy|$. Because we are
concerned with up to octopole moments of the
Green's function, we derive expressions for
first, second, and third order directional
derivatives of $G$ in the multipole basis
$\{Q_n,K_n\}$. As in \cref{sec:fundamental-solutions},
we make use of the standard suffix notation
and set $\v r = \xx - \yy$.

The following identities are straightforward
to derive. We have
\begin{align}
  \di G(\v r) &= G'(r) \frac{r_i}{r}, \\
  \di \dj G(\v r) &= \frac{G'(r)}{r} \delta_{ij} 
                    + \left( G''(r) 
                    - \frac{G'(r)}{r} \right) 
                    \frac{r_i r_j}{r^2},\\
  \begin{split}
    \di \dj \dk G(r) &= \left(\frac{G''(r)}{r} -
      \frac{G'(r)}{r^2}\right) \frac{ \delta_{jk} r_i + \delta_{ik}
      r_j
      + \delta_{ij} r_k }{r} \\
    & \quad + \left(G'''(r) - 3 \frac{G''(r)}{r} + 3
      \frac{G'(r)}{r^2}\right)
    \frac{r_i r_j r_k}{r^3} \; .
  \end{split}\\
\end{align}
From the recurrence and derivative relations for
modified Bessel functions \cite[S10.29]{NIST:DLMF},
we obtain

\begin{align}
  G'(r) &= -\frac{Q_1(r)}{2\pi\alpha} \; ,\\
  \frac{G'(r)}{r} &= \frac{Q_2(r)-K_0(\alpha r)}{4\pi} \; ,\\
  G''(r) - \frac{G'(r)}{r} &= -\frac{Q_2(r)}{2\pi} \; ,\\
  \frac{G''(r)}{r} - \frac{G'(r)}{r^2} &=
 -\alpha \frac{Q_3(r) - K_1(\alpha r)}{8\pi} \, ,\\
 G'''(r) - 3 \frac{G''(r)}{r} + 3 \frac{G'(r)}{r^2}
 &= \alpha \frac{Q_3(r)}{2\pi} \; .
\end{align}

To derive the separation of variables form
of the directional derivatives, we identify the
point $\v r$ with $r e^{i\theta}$. In the following,
we switch freely between identifying vectors by
their suffix and complex valued notations.
Let $\v \lambda, \v \mu, \v \nu$ be vectors
and $\lambda, \mu, \nu$ denote their complex
valued counterpart. Then,

\begin{align}
  \frac{r_i \lambda_i}{r} &= \frac{e^{i\theta}\bar\lambda
    + e^{-i\theta}\lambda}{2} \; ,
\end{align}
so that

\begin{equation}
  \label{eq:dipole}
  \partial_\lambda G = -\frac{Q_1(r)}{4\pi\alpha}\left
  (e^{i\theta}\bar\lambda + e^{-i\theta}\lambda \right ) \; .
\end{equation}
Similarly, we have that

\begin{align}
  \frac{r_ir_j\lambda_i\mu_j}{r^2} &=
  \frac{e^{i2\theta}\bar \lambda \bar \mu + e^{-i2\theta} \lambda \mu
    + \lambda \bar \mu + \bar \lambda \mu }{4} \\
  \delta_{ij} \lambda_i \mu_j &= \frac{\lambda \bar \mu + \bar \lambda \mu}{2}  \; ,
\end{align}
so that 

\begin{equation}
  \label{eq:quadrupole}
  \partial_\lambda \partial_\mu G =
  -\frac{Q_2(r)}{8\pi}\left
  (e^{i2\theta}\bar\lambda\bar\mu + e^{-i2\theta}\lambda\mu \right )
  - \frac{K_0(\alpha r)}{8\pi} \left (
  \lambda \bar \mu + \bar \lambda \mu \right )\; .
\end{equation}
Finally, we have that

\begin{align}
  \frac{r_ir_jr_k\lambda_i\mu_j\nu_k}{r^3} &=
  \frac{e^{i3\theta}\bar\lambda\bar\mu\bar\nu +
    e^{-i3\theta}\lambda\mu\nu +
      e^{i\theta}(\bar\lambda\bar\mu\nu + \bar\lambda\mu\bar\nu
      + \lambda\bar\mu\bar\nu) +
      e^{-i\theta}(\lambda\mu\bar\nu + \lambda\bar\mu\nu +
      \bar\lambda\mu\nu)}{8} \\
  \frac{\delta_{jk}r_i + \delta_{ki}r_j + \delta_{ij}r_k}{r}
  \lambda_i\mu_j\nu_k &= \frac{e^{i\theta} (\bar\lambda\bar\mu\nu
    + \bar\lambda\mu\bar\nu + \lambda\bar\mu\bar\nu) +
    e^{-i\theta}(\lambda\mu\bar\nu+\lambda\bar\mu\nu
    +\bar\lambda\mu\nu)}{2} \; ,
\end{align}
so that 

\begin{equation}
  \label{eq:octopole}
  \partial_\lambda \partial_\mu \partial_\nu G =
  \frac{\alpha Q_3(r)}{16\pi}\left
  (e^{i3\theta}\bar\lambda\bar\mu\bar\nu +
  e^{-i3\theta}\lambda\mu\nu \right )
  + \frac{\alpha K_1(\alpha r)}{16\pi} \left (
e^{i\theta} (\bar\lambda\bar\mu\nu
    + \bar\lambda\mu\bar\nu + \lambda\bar\mu\bar\nu) +
    e^{-i\theta}(\lambda\mu\bar\nu+\lambda\bar\mu\nu
    +\bar\lambda\mu\nu) \right )  \; .
\end{equation}
\begin{remark}
  Note that the $e^{in\theta}$ terms always have
  a coefficient which is the complex conjugate
  of the $e^{-in\theta}$ terms. This is used by the
  fast multipole method to reduce the storage for
  expansions and to reduce the computational
  cost of translating expansions.
\end{remark}

Once the exact separation of variables representation
for a charge or higher order moment of $G$ is obtained
from the formulas above, one can simply use the
translation formulas for multipole expansions to shift
each of the per-source expansions to the center of the
leaf box. The cost of this operation depends
on the number of sources, but is only performed
on leaf boxes. Similarly, the per-source
expansions can be translated, transformed,
and added to the local expansion in any box
for which the sources are well-separated.

%%%%%%%%%%%%%%%%%%%%%%%%%%%%%%%%%%%%%%%%%%%%%%%%%%%%%%%%%%%%%%%%%%%%%%%%%%%%%
\subsubsection{Fast direct solver}
\label{sec:fast-direct-solver}

Let $\Amat$ be the $2\Nbdry \times 2\Nbdry$ matrix
corresponding to the discretized boundary
integral equation, \cref{eq:BIE_quad}, so that

\begin{equation}
  A_{2(i-1)+k,2(j-1)+l} =
  \frac{1}{2} \delta_{2(i-1)+k,2(j-1)+l}
  + D_{kl}(\xx_i,\xx_j) w_{j}(\xx_i) \; ,
\end{equation}
with $1 \leq i,j \leq \Nbdry$ and $1\leq k,l \leq 2$.
Na\"ive inversion of this matrix would require
$O(\Nbdry^3)$ floating point operations. However,
the complexity of inversion can be reduced
to $O(\Nbdry)$ by taking advantage of structure within
the matrix.

The ability to use a truncated multipole expansion
to represent well-separated interactions to
high precision, as in the previous section,
implies that sub-matrices of $\Amat$ for which the
row points are well-separated from the column
points are of low numerical rank (and likewise
when the column points are well-separated from the
row points).
Because of these low rank structures in the
matrix, there is an ordering of the indices for
which the matrix $\Amat$ is hierarchically
block separable (HBS), as in \cite{Gillman2012,Marple2016}.
For such matrices, it is possible to efficiently
compute a compressed representation of the matrix
in $O(\Nbdry)$ operations and using $O(\Nbdry)$
storage. Then, a representation of the inverse
can also be computed and applied in $O(\Nbdry)$
operations based on that compressed representation.
Such methods combine matrix compression \cite{cheng_2005}
with a variant of the Sherman-Woodbury-Morrison
formula and take advantage of hierarchical
features of the matrix, as in the FMM.
See \cite{Gillman2012,Marple2016} for more
details.

A common feature of fast direct solvers is that
there is a relatively large upfront cost
in forming the representation of the inverse
while the cost of the subsequent linear system
solves for any given right hand side is orders of
magnitude smaller. This is well-suited to the
numerical examples in this paper, where
we have a fixed boundary discretization and
use a fixed time step length. We only pay the upfront
cost once and use the same representation of the
inverse for fast system solves at each step in the
simulation.

\begin{remark}
  Technically, the $O(\Nbdry)$ scaling of the
  HBS method depends on a
  more subtle aspect of the low-rank
  structure of the matrix.
  Because the boundary
  points are contained to a one dimensional curve,
  even the interactions between adjacent sections
  of the curve are of relatively low rank.
  This is illustrated in \cref{fig:lowrank_adjacent_1d}
  for two adjacent sections of the boundary, one
  gray and one black. Observe that the points on the
  gray portion of curve are approximately well-separated
  from each set of points between the dotted separators
  on the black portion of the curve.
%  Note that by dividing the black portion of the
%  curve in half, then cutting the left half in
  %  half again, and so on the points on the
%  gray portion of curve are approximately well-separated
%  from each set of points between the dotted separators
%  on the black portion of the curve. If we end the
%  cutting-in-half procedure when there is some
%  small number of discretization nodes remaining
%  in the left-over section of the black curve, then
%  we see that, assuming the points are roughly evenly
%  spaced, the number of required cuts grows only
%  logarithmically in the number of points on the black
%  section of the curve. Because the rank of interaction
%  for each section is low and there are relatively
%  few of them, the overall rank of interaction is low
%  for these adjacent pieces of curve.
\end{remark}

%%%%%%%%%%%%%%%%%%%%%%%%%%%%%%%%%%%%%%%%%%%%%%%%%%%%%%%%%%%%%%%%%%%%%%%%%%%%%
\subsection{Time stepping} 

The integral equation approach that we have outlined so far is based
on a first order IMEX discretization of the Navier-Stokes equations
\eqref{eq:navier-stokes-imex}. In order to increase the temporal order
of accuracy, we use the semi-implicit spectral deferred corrections
(SISDC) method \cite{Minion2003}. Consider an autonomous initial value
problem of the form
\begin{align}
  \phi'(t) &= G(\phi(t)) = G_E(\phi(t)) + G_I(\phi(t)), \\
  \phi(0) &= \phi_0,
\end{align}
where $G_I$ is a stiff part that needs implicit treatment, and $G_E$
is a nonlinear term that would be impractical to treat implicitly. The method
is based on dividing a time interval $[t_a, t_b]$,
$t_b - t_a = \Delta t$, into $p$ subintervals by choosing points $t_m$
such that $t_a = t_0 < t_1 < \dots < t_p = t_b$, with
$\deltat_m = t_{m+1} - t_m$. A first approximation of $\phi(t)$ on
$[t_a, t_b]$, denoted $\phi^0(t)$, is computed at the substeps
$\phi_m^0 = \phi^0(t_m)$ using a first order IMEX discretization,
\begin{align}
  \phi_{m+1}^0 - \delta t_m G_I(\phi_{m+1}^0) &= \phi_m^0 + \deltat_m G_E(\phi_m^0),
                                                \quad m=0,\dots,p-1,
                                                \label{eq:sisdc_initial}\\
  \phi_0^0 &= \phi(t_a).
\end{align}
Then, a sequence of increasingly accurate approximations
$\{\phi^1(t), \phi^2(t), \phi^3(t), \dots\}$ are computed at the nodes
$\phi_m^k=\phi^k(t_m)$ using the correction equation
\begin{align}
  \begin{split}
    \phi_{m+1}^{k+1} - \delta t_m G_I(\phi_{m+1}^{k+1}) = \phi_m^{k+1} + \deltat_m \left[
      G_E(\phi_m^{k+1}) - G_E(\phi_m^{k})
       - G_I(\phi_{m+1}^{k}) \right]
     + I_m^{m+1}(\phi^k_0, \dots, \phi^k_p), \\
     m=0, \dots, p-1, \\
     k=0, 1, 2, \dots,     
   \end{split}
  \label{eq:sisdc_formula}  
\end{align}
all initialized to $\phi^k_0=\phi(t_a)$.  The computation of
$\{\phi_0^{k+1},\dots,\phi_p^{k+1}\}$ is in this context considered
one application of the correction equation. The term
$I_m^{m+1}(\phi^k_0, \dots, \phi^k_p)$ is a quadrature approximation
of the integral
\begin{align}
  I_m^{m+1}(\phi^k) = \int_{t_m}^{t_{m+1}} G(\phi^k(t)) \dif t. 
\end{align}
The quadrature is precomputed, so that computing the
integral requires only a small matrix-vector multiplication,
as in
\begin{align}
  I_m^{m+1}(\phi^k_0, \dots, \phi^k_p) &= \sum_{j=0}^p q_j^{m} G(\phi^k_j), 
                                         \quad m=0,\dots,p-1 .
\end{align}
The weights $q_j^m$ are chosen so that they integrate the
interpolating polynomial through the $p+1$ points, which gives an
$\mathcal O(\Delta t^{p+2})$ error in each quadrature. In
\cite{Minion2003} the points are chosen to be the nodes of
Gauss-Lobatto quadrature. This choice ensures that the polynomial
interpolation is well-conditioned for high order schemes, while
conveniently including nodes at the beginning and end of the
interval.

Each application of the correction equation increases the local order
of accuracy by one, so that if
$\phi(t_b) -\phi^0_p = \mathcal O(\Delta t^2)$, then after $k$
corrections we have
\begin{align}
  \phi(t_b) -\phi^k_{p} = \mathcal O( \Delta t^{k+2} + \Delta t^{p+2}) .
\end{align}
% number of substeps p = sdc_order - 1
% number of corrections K = sdc_order - 1
% solves for initial: p
% solves per correction: p
% => total work: p*(K+1) = (sdc_order-1)*sdc_order
A SISDC time stepping method with global order of accuracy $K$
therefore requires $p=K-1$ substeps and $K-1$ applications of the
correction equation. Taking into account the $p$ steps taken to
compute the first approximation \eqref{eq:sisdc_initial}, the total
number of implicit linear system solves needed to take one step of
length $\Delta t$ is then $K^2-K$, resulting in the final
approximation $\phi(t_b) \approx \phi_{K-1}^{K-1}$.

In our case, the equation of interest is the momentum equation of the
Navier-Stokes equations \eqref{eq:navier-stokes}. We split the
equation into implicit and explicit parts as before, setting
\begin{align}
  G_E &= -(\uu \cdot \nabla) \uu, \\
  G_I &= \frac{1}{\Re}\Delta \uu - \nabla p .
\end{align}

Our scheme for the modified Stokes equation relies on precomputing a
solver and quadrature rules for a specific $\alpha$, which in this
context is a function of the substep length. To reduce the
precomputation time and storage, it is therefore preferable for us to
use substeps that are of equal length $\deltat = \Delta t /
(K-1)$. This coincides with the points used by \cite{Minion2003} for
$K=2$ (where $\deltat=\Delta t$) and $K=3$ (since 3-point
Gauss-Lobatto is equidistant). For moderate values of $K>3$, equispaced
substeps could still be feasible, though at some level of $K$ the
ill-conditioning of polynomial interpolation on equispaced nodes would
make itself known.

Assuming a constant substep length of $\deltat$, the first
approximation is computed by using our composite scheme to solve the
modified Stokes equation \eqref{eq:modstokes},
% \begin{align}
%   \uu^0_{m+1} - \deltat 
%   \left[ \frac{1}{\Re}\Delta \uu^0_{m+1} - \nabla p^0_{m+1} \right] 
%   &=  \uu^0_m -\deltat (\uu^0_N \cdot \nabla) \uu^0_N,
% \end{align}
\begin{align}
  (\alpha^2  - \Delta) \uu^0_{m+1} + \Re \nabla p^0_{m+1} &=  \FF^0_{m+1}, \quad \alpha^2={\Re/\deltat},
\end{align}
with the forcing function from the original IMEX scheme
\eqref{eq:navier-stokes-imex},
\begin{align}
  \FF^0_{m+1} = \alpha^2 \uu^0_m - \Re (\uu^0_N \cdot \nabla) \uu^0_N .
  \label{eq:rhs_initial}
\end{align}
To apply the correction equation, we solve the same PDE,
\begin{align}
  (\alpha^2 - \Delta ) \uu^{k+1}_{m+1} + \Re \nabla p^{k+1}_{m+1} = \FF_{m+1}^{k+1},
\end{align}
this time with the right hand side
\begin{align}
  \FF_{m+1}^{k+1} = \alpha ^2 \uu_m^{k+1} - \Re \left[ (\uu_m^{k+1} \cdot \nabla)
                    \uu_m^{k+1} - (\uu_m^{k} \cdot \nabla) \uu_m^{k} \right] 
                    -\Delta \uu^{k}_{m+1} + \Re \nabla p^{k}_{m+1} +
                    \frac{\Re}{\deltat} I_m^{m+1}, 
                    \label{eq:rhs_correction}
\end{align}
and the quadrature term
\begin{align}
\frac{\Re}{\deltat}  I_m^{m+1} &= \frac{1}{\deltat} \sum_{j=0}^p q_j^{m} \left[
                                   \Delta \uu^k_j - \Re \nabla p^k_j - \Re (\uu^k_j \cdot \nabla) \uu^k_j
                                   \right] 
                                   \label{eq:rhs_correction_integral} .
\end{align}

For the first order IMEX scheme, we need to compute the solution $\uu$
and its gradient $\nabla\uu$ on the grid at every substep, and we also
need them at the initial time $t=0$ to start the simulation. In
addition to that, the SISDC scheme also requires the quantity
$\Delta \uu - \Re \nabla p$ at every substep, and at the initial
time. At a substep $(m,k)$ we can compute this quantity without
actually having to compute $\Delta\uu_m^k$ and $\nabla p_m^k$
separately, simply by recovering it from the equation solved to
compute the solution $\uu^k_m$ at that substep,
\begin{align}
  \Delta  \uu^{k}_{m} - \Re \nabla p^{k}_{m} = \alpha^2  \uu^{k}_{m} - \FF_{m}^{k} .
\end{align}
For the initial time, this quantity is not available. Note that,
when added to the right hand side, any irrotational term (which is the
gradient of a scalar function) has no effect on the computed velocity
field. Therefore, we only require the value of $\Delta \uu$ be supplied
at the initial time and use $\Delta \uu^{0}_0$ instead
of $\Delta \uu^{0}_0 - \Re \nabla p^0_0$ in correction
equations for the first interval.

\begin{remark}
  The fact that the irrotational part of the right hand side
  has no effect on the velocity is clear in two
  dimensions when analyzing the stream function formulation of
  the modified Stokes equations, \cref{eq:modstokes}.
  Let $\uu = \nabla^\perp \Psi$ and write
  the Helmholtz decomposition of the right hand side $\FF$ as
  $\FF = \nabla \Phi + \nabla^\perp \Xi$ for some scalar functions
  $\Psi$, $\Phi$, and $\Xi$. Substituting these expressions into
  \cref{eq:modstokes} and taking the dot product of both sides
  with $\nabla^\perp$, we obtain

  \[ (\alpha^2 - \Delta) \Delta \Psi = \Delta \Xi \; ,\]
  which depends only on the solenoidal part of $\FF$.
\end{remark}
  
%%% Local Variables: 
%%% TeX-master: "nseIE"
%%% End: 

%%%%%%%%%%%%%%%%%%%%%%%%%%%%%%%%%%%%%%%%%%%%%%%%%%%%%%%%%%%%%%%%%%%%%%%%
\section{Quadrature details}
\label{sec:quadrature}

As mentioned in \cref{sec:homogeneous-solution}, composite
Gauss-Legendre quadrature over $\bdry$ is only accurate for target
points that are well-separated from $\bdry$. More precisely, the
contribution to the double layer potential from a particular panel
$\Gamma$, using the direct Gauss-Legendre quadrature, is only accurate
if the target point $\xx$ is sufficiently far away from that
panel. Our goal, then, is to have a quadrature that can accurate
evaluate
\begin{align}
  \int_\panel \mu_i(\yy) T_{ijk}(\xx , \yy) \nhat_k(\v y)  \dSy,
\end{align}
for any panel $\Gamma\subset\bdry$ and target point $\xx$.  For this, we will
distinguish between two different cases: the singular case, when
$\xx \in \bdry$, and the nearly singular case, when $\xx \in \domain$
is close to $\bdry$.

Our method of choice for the singular and nearly singular quadrature
will be the kernel-split, complex interpolatory quadrature scheme of
Helsing et al. \cite{Helsing2008,Helsing2009,Helsing2015}. This scheme
is both accurate and efficient. In addition, at each target point it
only requires a finite (and small) number of operations to correct the
global direct quadrature returned by the FMM (it is
\emph{FMM-compatible}, according to the definition of \cite{Hao2014}).

\subsection{Kernel split}
\label{sec:kernel-split}

The core problem that kernel-split quadrature addresses is that 
the stresslet contains several singularities as $\xx\to\yy$. In
order to apply a kernel-split quadrature, we must first rewrite the
kernel in a way that explicitly exposes these singularities. The
closed-form expression for the stresslet, derived in
\cref{sec:fundamental-solutions}, is
\begin{align}
  \label{eq:stresslet_explicit}
  T_{ijk}(\v x, \v y) &= \alpha^2 \mathcal T_1(\alpha \norm{\rr}) \left( 
                  \delta_{jk} r_i 
                  + \delta_{ik} r_j 
                  + \delta_{ij} r_k \right)
                  + \alpha^4 \mathcal T_2(\alpha \norm{\rr}) 
                  r_i r_j r_k
                 + \alpha^2 \mathcal T_3(\alpha \norm{\rr}) \delta_{ik} r_j,
\end{align}
where $\rr = \xx - \yy$, and the stresslet functions
$\mathcal T_1$--$\mathcal T_3$ can be expressed in terms of modified
Bessel functions of the second kind,
\begin{align}
  \mathcal T_1(z) &= -\frac{ 2 z^2 K_0(z)+\left(z^2+4\right) z K_1(z)-4}{2 \pi  z^4}, \label{eq:stresslet_T1}\\
  \mathcal T_2(z) &= \frac{4 z^2 K_0(z)+\left(z^2+8\right) z K_1(z)-8}{\pi  z^6}, \label{eq:stresslet_T2}\\
  \mathcal T_3(r) &= \frac{z K_1(z)-1}{2 \pi  z^2}. \label{eq:stresslet_T3}
\end{align}
It is clear that any singular behavior in the stresslet must come from
the functions $\mathcal T_i$. Using a standard result
\cite[\S10.31]{NIST:DLMF}, we can decompose, or split, $K_0$ and $K_1$ into explicit
singularities as
\begin{align}
  K_0(z) &= K_0^S(z) - I_0(z) \log z, \\
  K_1(z) &= K_1^S(z) + I_1(z) \log z + \frac{1}{z}.
\end{align}
Here $K_0^S(z)$ and $K_1^S(z)$ are simply the remainders after
removing the singularities; both of these terms, along with $I_0$ and $I_1$
(modified Bessel functions of the first kind), are smooth functions of
$z$.  Using this, we can now split the stresslet functions
$\mathcal T_i$ into smooth functions multiplying explicit
singularities,
\begin{align}
  \mathcal T_1(z) &= \mathcal T_1^S (z) + \mathcal T_1^L (z) \log z, \label{eq:T1_split}\\
  \mathcal T_2(z) &= \mathcal T_2^S (z) + \mathcal T_2^L (z) \log z 
                    + \frac{1}{8 \pi z^2}
                    - \frac{1}{\pi z^4}, \label{eq:T2_split}\\
  \mathcal T_3(z) &= \mathcal T_3^S (z) + \mathcal T_3^L (z) \log z, \label{eq:T3_split}
\end{align}
where
\begin{align}
  \mathcal T_1^S(z) &= -\frac{2 z K_0^S(z)+\left(z^2+4\right) K_1^S( z)+z}{2 \pi  z^3}, &
  \mathcal T_1^L(z) &= \frac{2 z I_0(z)-\left(z^2+4\right) I_1(z)}{2 \pi  z^3}, 
  \\
  \mathcal T_2^S(z)  &= \frac{32 z K_0^S(z)+8 \left(z^2+8\right)
                       K_1^S(z) - z \left(z^2-16\right)}{8 \pi  z^5} &
  \mathcal T_2^L(z) &= \frac{\left(z^2+8\right) I_1(z)-4 z I_0(z)}
                      {\pi  z^5}, 
  \\
  \mathcal T_3^S(z) &= \frac{K_1^S(z)}{2 \pi z} &
  \mathcal T_3^L(z) &= \frac{I_1(z)}{2 \pi z}.
\end{align}
\begin{figure}
  \centering
  % Generated by: julia/output_scripts/kernel_split.jl
  \includegraphics[width=0.49\textwidth]{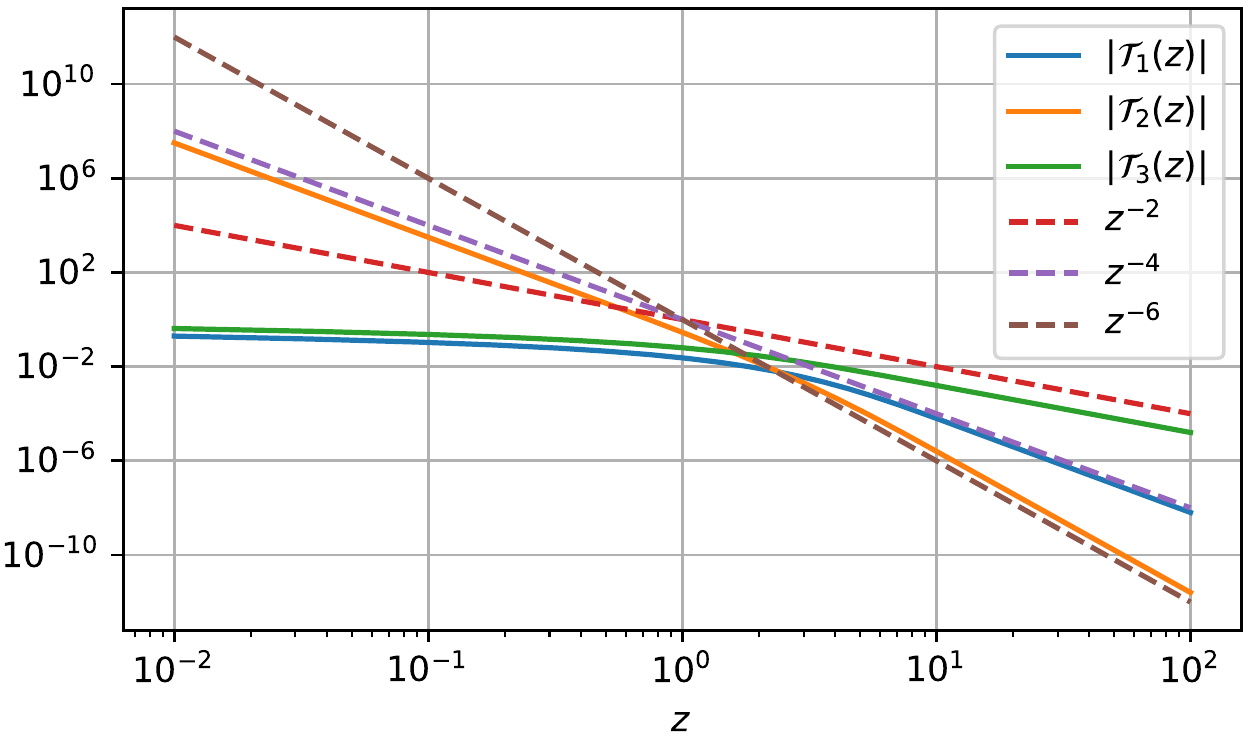}
  \includegraphics[width=0.49\textwidth]{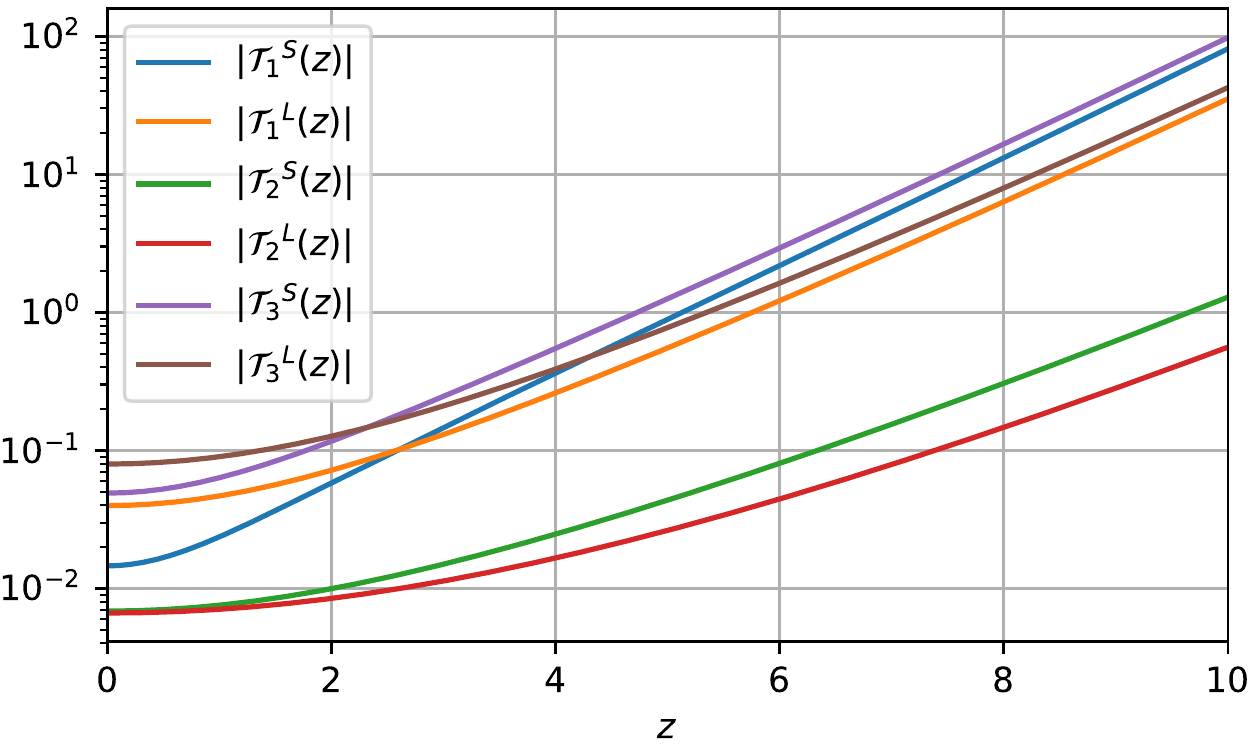}
  \caption{Left: The stresslet functions $\mathcal
    T_1$--$\mathcal T_3$, as defined in
    \cref{eq:stresslet_T1,eq:stresslet_T2,eq:stresslet_T3}. Right: The
    splits $\mathcal T_1^S$--$\mathcal T_3^L$, as defined in
    \cref{eq:T1_split,eq:T2_split,eq:T3_split}. Note that the
    magnitudes of the stresslet functions decay algebraically as
    $z\to\infty$, even though the magnitudes of their decompositions
    grow exponentially.}
  \label{fig:kernel_split}
\end{figure}
All of the above expressions are prone to cancellation for small
values of $z=\alpha\norm{\v r}$, and must then be evaluated using
power series expansions, discussed further in
\cref{sec:power-series}.
% From those expansions it can also be
% determined that all of the smooth functions $\mathcal T_i^S$ and
% $\mathcal T_i^L$ have constant limits as $z \to 0$.

Using the above definitions, we now split the stresslet as
\begin{align}
  T_{ijk}(\v r) = 
  \left( T_{ijk}^S(\v r) + T_{ijk}^L(\v r) \log(\alpha) \right)
  + T_{ijk}^L(\v r) \log\norm\rr
  + T_{ij}^C(\v r) \frac{r_k}{\norm\rr^2} 
  + T^Q \frac{r_i r_j r_k}{\norm\rr^4} ,
  \label{eq:stresslet_split}
\end{align}
where 
\begin{align}
  T_{ijk}^S(\v r) &= \alpha^2 \mathcal T_1^S(\alpha \norm\rr) 
                    \left( 
                    \delta_{jk} r_i 
                    + \delta_{ik} r_j 
                    + \delta_{ij} r_k \right)
                    + \alpha^4 
                    \mathcal T_2^S(\alpha \norm\rr) 
                    r_i r_j r_k
                    + \alpha^2 
                    \mathcal T_3^S(\alpha \norm\rr) \delta_{ik} r_j,
                    \label{eq:TS}
  \\
  T_{ijk}^L(\v r) &= \alpha^2 \mathcal T_1^L(\alpha \norm\rr) 
                  \left( 
                  \delta_{jk} r_i 
                  + \delta_{ik} r_j 
                  + \delta_{ij} r_k \right)
                  + \alpha^4 \mathcal T_2^L(\alpha \norm\rr) 
                  r_i r_j r_k
                    + \alpha^2 \mathcal T_3^L(\alpha \norm\rr) \delta_{ik} r_j,
                    \label{eq:TL}
  \\
  T_{ij}^C(\v r) &= \alpha^2 r_i r_j / 8\pi, \label{eq:TC}\\
  T^Q &= -1 / \pi, \label{eq:TQ}
\end{align}
Using these results, together with the relation
$T(\xx,\yy)=-T(\yy,\xx)$, we write the double layer potential
contribution \eqref{eq:dbl_lyr_pot} from a panel $\Gamma$ as
\begin{align}
  \int_\panel \mu_i(\yy) \left(
  T_{ijk}^S(\v r) + T_{ijk}^L(\v r) \log\norm\rr
  + T_{ij}^C(\v r) \frac{r_k}{\norm\rr^2} 
  + T^Q \frac{r_i r_j r_k}{\norm\rr^4}                            
  \right)\nhat_k(\v y)  \dSy, \quad \v r = \yy - \xx .
  \label{eq:dbl_lyr_pot_split}
\end{align}
This form, where the integrand is a smooth function plus a sum of
smooth functions multiplied by explicit singularities, is what allows
us to design efficient quadrature.

\subsection{On-boundary evaluation of the double layer potential}
\label{sec:boundary-evaluation}

When solving the discretized integral equation \eqref{eq:BIE_quad}, we
need to evaluate the double layer potential for target points $\xx$ on
the boundary $\bdry$. If $\xx$ is close to $\Gamma$ but distant in
arclength, then $\xx$ is treated as a nearly singular point, which we
will return to in \cref{sec:near-bound-eval}. This can happen,
for example, if $\bdry$ is multiply-connected or it has a narrow
channel. If on the other hand $\xx$ is close to $\Gamma$ both in
distance and arclength, which in practice means that $\xx$ is on
$\Gamma$ or on a neighboring panel, then we must consider what happens
with the singularities in the double layer potential
\eqref{eq:dbl_lyr_pot_split} in the limit $\xx\to\yy$ on $\bdry$: the
$\log\norm\rr$-term is singular in the limit, but that singularity is
cancelled by $T^L\to 0$. The $\rr\cdot\v\nhat/\norm\rr^2$-term has a
finite limit, which is then removed by $T^C\to 0$. The smooth term
$T^S$ also disappears, so all that remains is the smooth limit of the
fourth term,
\begin{align}
  \lim_{\substack{\v x \to \v y\\\v x, \v y \in \bdry}} 
  T_{ijk}(\v y - \v x) \nhat_k(\v y) = 
  \frac{\that_i(\v y) \that_j(\v y) \kappa(\v y)}{2 \pi},
\end{align}
where $\v\that=(\nhat_2, -\nhat_1)$ is the unit tangent vector and
$\kappa$ is the local curvature. However, even though the
log-singularity is cancelled out by its multiplying factor, we must
still apply kernel-split quadrature to it, in order to evaluate it to
high accuracy. This is done in precisely the same way as for
near-boundary evaluation in \cref{sec:near-bound-eval}. The remaining
three terms in \eqref{eq:dbl_lyr_pot_split} can be evaluated directly
using Gauss-Legendre quadrature, as they are smooth on $\bdry$.

\subsection{Near-boundary evaluation of the double layer potential}
\label{sec:near-bound-eval}

For a target point $\xx\in\domain$ close to $\Gamma$, the kernel of
\eqref{eq:dbl_lyr_pot_split} has singularities close to the interval
of integration, which are often referred to as near-singularities in the
literature. This makes the integrand a bad
candidate for polynomial approximation, which is why direct
Gauss-Legendre quadrature fails. The kernel-split quadrature, by
contrast, computes weights for each term of
\eqref{eq:dbl_lyr_pot_split} separately, and works well as long as the
function multiplying the singularity is smooth. The scheme works by
first rewriting the singularities in complex form, resulting in
linear combinations of the integrals
\begin{align}
  \int_\panel f(\tau) \log(\tau-z) \dif \tau 
  \quad \text{ and } \quad 
  \int_\panel f(\tau) \frac{ \dif\tau }{ (\tau-z)^m }, \: m=1,2,\dots.
  \label{eq:complex_forms}
\end{align}
Here $z = x_1 + i x_2$ and $\tau = y_1 + i y_2$ are complex
representations of the target point and the panel, and $f$ is assumed
to be a smooth, real-valued function. Then, for a panel with $n$
nodes, $f$ is approximated using a monomial expansion,
\begin{align}
  f(\tau) \approx \sum_{k=1}^{n} c_k \tau^{k-1},
\end{align}
such that, e.g., 
\begin{align}
  \int_\panel f(\tau) \frac{ \dif\tau }{ (\tau-z)^m } \approx
  \sum_{k=1}^{n} c_k p_k,
  \quad \text{ where } \quad
  p_k = \int_\panel \frac{ \tau^{k-1} \dif\tau }{ (\tau-z)^m } .
  \label{eq:specquad_pk}
\end{align}
The integrals $p_k$ can be evaluated using recursion formulas. By
solving a transposed $n\times n$ Vandermonde system, quadrature
weights $w_j$ can be computed \cite{Helsing2008}, such that the integral
of a function $f$ known at the nodes $\tau_j$ can be approximated as
\begin{align}
  \int_\panel f(\tau) \frac{ \dif\tau }{ (\tau-z)^m } \approx
  \sum_{j=1}^n w_j f(\tau_j) .
\end{align}
This quadrature is $n$th order accurate, as it is based on
interpolation of $f$ at the $n$ quadrature nodes. Solution of the
complex $n \times n$ Vandermonde system can be accelerated using the
Bj\"orck-Pereyra algorithm \cite{Bjorck1970}.

In our layer potential \eqref{eq:dbl_lyr_pot_split}, the first term is
smooth and can be evaluated using the direct Gauss-Legendre
quadrature. The remaining three terms must be converted into the forms
\eqref{eq:complex_forms} before we can evaluate them. For this, we use
the following notation:
\begin{align*}
    z &= x_1 + i x_2, &
    \tau(t) &= y_1(t) + i y_2(t), \\
    \nu_\tau &= i \tau'(t) / \abs{\tau'(t)}, &
    \omega &= \mu_1 + i\mu_2 .
\end{align*}
Here $\yy(t)$ is the parametrization of $\Gamma$, as defined in
\ref{sec:comp-homog-corr}. Identifying $\mathbb R^2$ with
$\mathbb C$ and using $\rr=\yy-\xx=\tau-z$, the remaining three terms
can now be cast into complex forms using the relations
\begin{align}
  \int_\panel f(\xx,\yy) \log\norm{\v r} \dif S
  &= -\imag \int_\panel f(z,\tau) \conj\nu_\tau \log(\tau - z) \dif\tau, \label{eq:cpx_log}\\
  \int_\panel f(\xx,\yy) \frac{ \v r \cdot \v\nhat }{\norm{r}^2} \dif S
  &= -\imag \int_\panel f(z,\tau) \frac{ \dif \tau }{\tau - z}, \label{eq:cpx_cauchy}\\
  \int_\Gamma \frac{ (\v \mu \cdot \v r) \v r (\v r \cdot \v\nhat)}{\norm{\v r}^2} \dif S 
  &= 
    -\frac{1}{2}\int_\Gamma \omega\imag[(\conj\tau - \conj z) \dif \tau]
    -\frac{1}{4i}\int_\Gamma \conj\omega(\tau-z) \dif \tau 
    +\frac{1}{4i} \overline{\int_\Gamma \frac{ \omega(\conj \tau-\conj z)^2}{(\tau - z)}\dif\tau},
  \label{eq:cpx_cauchy2} \\
  \int_\Gamma \frac{ (\v \mu \cdot \v r) \v r (\v r \cdot \v\nhat)}{\norm{\v r}^4} \dif S 
  &= -\frac{1}{2}\int_\Gamma \omega \imag\left[ \frac{ \dif \tau}{(\tau - z} \right] 
    +\frac{1}{4i}\overline{\int_\Gamma \frac{\omega \bar\nu^2 \dif \tau}{\tau - z}}
    +\frac{1}{4i}\overline{ \int_\Gamma \frac{\omega(\bar \tau-\bar z) \dif \tau }{(\tau - z)^2} }.
    \label{eq:cpx_stresslet}
\end{align}
In order to use the relation \eqref{eq:cpx_log} to compute the second
term of \eqref{eq:dbl_lyr_pot_split}, we set
$f=\mu_i T^L_{ijk} \nhat_k$ for $j=1,2$. The third term of
\eqref{eq:dbl_lyr_pot_split} can be computed either using
\eqref{eq:cpx_cauchy} with $f=\alpha^2\mu_i r_i r_j / 8\pi$, $j=1,2$,
or \eqref{eq:cpx_cauchy2}. We use \eqref{eq:cpx_cauchy}, since it
corresponds to the Laplace double layer potential, for which there is
code readily available in \cite{Helsing2015}. The fourth term of
\eqref{eq:dbl_lyr_pot_split} corresponds directly to
\eqref{eq:cpx_stresslet}, and is actually the stresslet of regular
Stokes flow, for which quadrature was developed in \cite{Ojala2015}.

By applying the above relations, and computing quadrature weights for
the kernels $\log(\tau-z)$, $(\tau-z)^{-1}$, and $(\tau-z)^{-2}$, we
can compute quadrature weights that accurately evaluate all the terms
of \eqref{eq:dbl_lyr_pot_split} for $\xx$ arbitrarily close to
$\Gamma$.

% Load gradient derivations from separate source
\subsection{Near-boundary evaluation of the gradient of the double layer potential}
\label{sec:gradient}

In order to get accurate values of the gradient of the double layer
potential, we need to apply the above procedure to a kernel-split form
of the integral
\begin{align}
  \dpd{}{x_l}
  \int_\panel \mu_i(\yy) \left(
  T_{ijk}^S(\v r) + T_{ijk}^L(\v r) \log\norm\rr
  + \frac{\alpha^2  }{ 8\pi } \frac{r_i r_j r_k}{\norm\rr^2} 
   -\frac{ 1 }{ \pi } \frac{r_i r_j r_k}{\norm\rr^4}                            
  \right)\nhat_k(\v y)  \dSy, \quad \v r = \yy - \xx .
  \label{eq:gradient_split_begin}
\end{align}

Direct differentiation of the first term of \eqref{eq:gradient_split_begin} gives
\begin{align}
  \begin{split}
    -\dpd{}{x_\ell} T_{ijk}^S(\v r) =& \frac{\alpha r_\ell}{\norm\rr}\left(
      \alpha^2 {\mathcal T}_1^{S,\prime}(\alpha \norm\rr) \left(
        \delta_{jk} r_i + \delta_{ik} r_j + \delta_{ij} r_k \right) +
      \alpha^4 {\mathcal T}_2^{S,\prime}(\alpha \norm\rr) r_i r_j r_k +
      \alpha^2 {\mathcal T}_3^{S,\prime}(\alpha \norm\rr) \delta_{ik}
      r_j \right) \\ &+ \alpha^2 {\mathcal T}_1^S(\alpha \norm\rr)
    \left( \delta_{jk} \delta_{i\ell} + \delta_{ik} \delta_{j\ell} +
      \delta_{ij} \delta_{k\ell} \right) \\ &+ \alpha^4
    {\mathcal T}_2^S(\alpha \norm\rr) \left( r_j r_k \delta_{i\ell} +
      r_i r_k \delta_{j\ell} + r_i r_j \delta_{k\ell} \right) +
    \alpha^2 {\mathcal T}_3^S(\alpha \norm\rr) \delta_{ik}
    \delta_{j\ell} .
    \label{eq:TS_deriv}
  \end{split}
\end{align}
This term is smooth, in spite of the leading $\norm\rr^{-1}$, since
the leading order term of ${\mathcal T}_n^{S,\prime}(\alpha \norm\rr)$
is $\mathcal O(\norm\rr)$ (see \cref{sec:power-series}). It can
therefore be evaluated using direct Gauss-Legendre quadrature.

Direct differentiation of the second term of \eqref{eq:gradient_split_begin} gives
\begin{align}
  \dpd{}{x_\ell} T_{ijk}^L(\v r) \log\norm\rr = 
  \dpd{ T_{ijk}^L(\v r) }{x_\ell} \log\norm\rr - T_{ijk}^L(\v r) \frac{r_l}{\norm\rr^2} .
  \label{eq:TL_log_deriv}
\end{align}
The derivative of $T^L$ is computed analogously to
\eqref{eq:TS_deriv}, and the first term of \eqref{eq:TL_log_deriv} is
then integrated using the same rule as for the $\log$ term of
\eqref{eq:dbl_lyr_pot_split}. The second term of
\eqref{eq:TL_log_deriv} is integrated using the complex interpolatory
quadrature for $(\tau-z)^{-1}$, after the rewrite
\begin{align}
  \int_\Gamma f(\xx,\yy) \frac{\yy-\xx}{\norm{\yy-\xx}^2} \dSy 
  = \overline{ \int_\Gamma  \frac{i f(z,\tau)}{\nu_\tau} \frac{\dif\tau}{\tau-z} },
\end{align}
using $f=\mu_i T^L_{ijk} \nhat_k$ for $j=1,2$.

In order to rewrite the third and fourth terms of
\eqref{eq:gradient_split_begin}, we will differentiate the complex
forms \eqref{eq:cpx_cauchy2} and \eqref{eq:cpx_stresslet}. In general,
if the velocity field is represented by a complex field $\phi$,
\begin{align}
  u_1 + i u_2 = \phi(z, \conj z),
\end{align}
then we can get the partial derivatives as
\begin{align}
  \dpd{}{x_1} (u_1 + i u_2) &= \dpd{\phi}{z} + \dpd{\phi}{\conj z}, &
  \dpd{}{x_2} (u_1 + i u_2) &= i \left( \dpd{\phi}{z} - \dpd{\phi}{\conj z} \right),
\end{align}
since
\begin{align}
  \dpd{z}{x_1} &= 
  \dpd{\conj z}{x_1} = 1, &
  \dpd{z}{x_2} &= i, &
  \dpd{\conj z}{x_2} &= -i .
\end{align}
Beginning with term three of \eqref{eq:gradient_split_begin}, we identify
\begin{align}
  \phi = 4 i \int_\Gamma \frac{\mu_i r_i r_j r_k \nhat_k}{\norm\rr^2} \dif S .
\end{align}
In complex form, using \eqref{eq:cpx_cauchy2}, we have
\begin{align}
  \phi(z, \conj z)
  &= 
    -\int_\Gamma \omega (\conj\tau - \conj z) \dif \tau
    + \int_\Gamma \omega (\tau - z) \dif\conj\tau
    -\int_\Gamma \conj\omega(\tau-z) \dif \tau 
    + \overline{\int_\Gamma \frac{ \omega(\conj \tau-\conj z)^2}{(\tau - z)}\dif\tau},
\end{align}
with partial derivatives that are straightforward to evaluate,
\begin{align}
  \dpd{\phi}{z}
  &=
    2i \imag \int_\Gamma \conj\omega \dif\tau
    -2\overline{\int_\Gamma \frac{ \omega(\conj \tau-\conj z)}{(\tau - z)}\dif\tau}
               ,\\
  \dpd{\phi}{\conj z}
  &=
    \int_\Gamma \omega \dif \tau
    + \overline{\int_\Gamma \frac{ \omega(\conj \tau-\conj z)^2}{(\tau - z)^2}\dif\tau}
    .
\end{align}
% ============== NOT USING PARTIAL INTEGRATION FOR THIS TERM
% We can reduce the singularity in the second term of $\pd{f}{\conj z}$ by partial integration:
% \begin{align}
%   \dpd{f}{\conj z}
%   &=
%     \int_\Gamma \omega \dif \tau
%     + \overline{ \left[
%     -\frac{ \omega(\conj \tau-\conj z)^2}{\tau - z}
%     \right]_{\tau_a}^{\tau_b} }
%     + \overline{\int_\Gamma \frac{
%     \od{\omega}{\tau}(\conj \tau-\conj z)^2 + 2\omega(\conj \tau-\conj z) \od{\conj\tau}{\tau}
%     }{\tau - z}\dif\tau}
% \end{align}

Lastly, for the fourth term of \eqref{eq:gradient_split_begin}, we identify
\begin{align}
  \phi = 4 i \int_\Gamma \frac{\mu_i r_i r_j r_k \nhat_k}{\norm\rr^4} \dif S .  
\end{align}
Written in complex form, using \eqref{eq:cpx_stresslet}, we have
\begin{align}
  4\pi i \phi(z, \conj z)
  &= 
    -\int_\Gamma \omega \frac{ \dif \tau}{\tau - z} 
    + \int_\Gamma \omega \frac{ \dif \conj\tau }{\conj{\tau} - \conj{z}} 
    +\overline{\int_\Gamma \frac{\omega \bar\nu^2 \dif \tau}{\tau - z}}
    +\overline{ \int_\Gamma \frac{\omega(\bar \tau-\bar z) \dif \tau }{(\tau - z)^2} } \\
  &= 
    -\int_\Gamma \frac{ \omega \dif \tau}{\tau - z} 
    +\overline{\int_\Gamma \frac{(\conj\omega + \omega \bar\nu^2) \dif \tau}{\tau - z}}
    +\overline{ \int_\Gamma \frac{\omega(\bar \tau-\bar z) \dif \tau }{(\tau - z)^2} } .
\end{align}
and
\begin{align}
  \dpd{\phi}{z} &= -\int_\Gamma \frac{ \omega \dif \tau}{(\tau - z)^2} 
               -\overline{ \int_\Gamma \frac{\omega \dif \tau }{(\tau - z)^2} } 
               = -2 \real \int_\Gamma \frac{ \omega \dif \tau}{(\tau - z)^2}  , \label{eq:dstressletdz} \\
  \dpd{\phi}{\conj z} &=
                     \overline{ \int_\Gamma \frac{(\conj\omega + \omega \bar\nu^2) \dif \tau}{(\tau - z)^2} }
                     + 2\overline{ \int_\Gamma \frac{\omega(\bar \tau-\bar z) \dif \tau }{(\tau - z)^3} } .
                        \label{eq:dstressletdzbar}
\end{align}
The last term has a singularity of type $(\tau-z)^{-3}$, which is one
order higher than what we have had to evaluate so far. While the
complex interpolatory quadrature scheme can in principle by used for
any order singularity, the accuracy can be expected to decrease with
higher order singularities. We therefore reduce the order of the
singularity by one through integration by parts, such that
\begin{align}
  \dpd{\phi}{\conj z} &=
                     -{
                     \left[
                        \frac{\conj\omega( \tau- z) }{(\conj\tau - \conj z)^2}
                     \right]_{\tau_a}^{\tau_b}
                     }
                     + \overline{ \int_\Gamma \frac{
                     \od{\omega}{\tau}(\bar \tau-\bar z) + \omega \od{\conj\tau}{\tau}
                     }{(\tau - z)^2}\dif \tau} .
                        \label{eq:dstressletdz_pi}
\end{align}
To evaluate this, we make use of the relation
$\od{\conj\tau}{\tau} = -\conj\nu^2$. The derivative
$\od{\omega}{\tau}$ is evaluated using polynomial interpolation and
differentiation on the panel.

% ============== FULL PARTIAL INTEGRATION
% Now we can use partial integration:
% \begin{align}
%   \dpd{\phi}{z} &= 2\real \left[\frac{\omega}{(\tau - z)}\right]_{\tau_a}^{\tau_b}
%                -2 \real \int_\Gamma \frac{ \od{\omega}{\tau} }{(\tau - z)} \dif \tau \\
%   \dpd{\phi}{\conj z} &=
%                      -\overline{
%                      \left[
%                      \frac{(\conj\omega + \omega \bar\nu^2)}{\tau - z}
%                      + \frac{\omega(\bar \tau-\bar z) \dif \tau }{(\tau - z)^2}
%                      \right]_{\tau_a}^{\tau_b}
%                      }
%                      +
%                      \overline{  \int_\Gamma \frac{                     
%                      \od{\conj\omega}{\tau} + \od{\omega}{\tau}\conj\nu^2 + \omega \od{\conj\nu^2}{\tau} 
%                      }{\tau - z} \dif \tau }
%                      + \overline{ \int_\Gamma \frac{
%                      \od{\omega}{\tau}(\bar \tau-\bar z) + \omega \od{\conj\tau}{\tau}
%                      }{(\tau - z)^2}\dif \tau} .
% \end{align}
% Some of the above derivatives can be recovered from the discretization by noting that
% \begin{align}
%   \dod{\conj\tau}{\tau} &= -\conj\nu^2, &
%   \dod{\conj\nu^2}{\tau} &= - \dod[2]{\conj\tau}{\tau} = 2 \frac{\kappa}{\nu^3} .
% \end{align}

%%% Local Variables: 
%%% TeX-master: "nseIE"
%%% End: 

\subsection{Implementation details}
\label{sec:impl-deta}

The quadrature described above is capable of evaluating the double
layer potential to high accuracy everywhere, both on and off the
boundary. However, the implementation must be carried out carefully in
order to obtain maximum accuracy and efficiency. We outline
the full scheme in this subsection.

\subsubsection{Determining where to apply kernel-split quadrature}
\label{sec:determ-where-apply}

A key element in implementing the kernel-split quadrature is
determining where to apply it, which translates into determining which
points are too close for the direct Gauss-Legendre quadrature, both on
and off the boundary. This is important, since setting the limit too
close to the source panel results in quadrature errors from the
singularities, while setting the limit too far away can cause
numerical errors in the kernel-split quadrature (although there is a
forgiving middle ground). To make this choice, we use the quadrature
error results in \cite{AfKlinteberg2018}. For a given source panel
$\Gamma$, we consider the mapping $\gamma : \mathbb{R} \to \mathbb{C}$
which takes the standard interval $[-1,1]$ to $\Gamma$. For integrals
of the type \eqref{eq:complex_forms}, the quadrature error at a target
point $z$ close to $\Gamma$ is dependent on the location of the
preimage of $z$ under a complexification of $\gamma$. For
$t^* \in \mathbb{C}$ such that $z = \gamma(t^*)$, the quadrature error
at $z$ is approximately proportional to $\rho(t^*)^{-2n}$, where
$\rho(t)=\abs{t \pm \sqrt{t^2-1}}$, with the sign defined such that
$\rho>1$, is the Bernstein radius. Even though it is possible to
derive more accurate quadrature estimates for the singularities of our
kernel \eqref{eq:dbl_lyr_pot_split} (see \cite{AfKlinteberg2016quad,
  AfKlinteberg2018}), it is for our purposes sufficient to
heuristically determine a limiting radius $R$, and then apply kernel
split quadrature for points $z$ such that $\rho(\gamma^{-1}(z)) <
R$.
Note that finding $t^* = \gamma^{-1}(z)$ is a cheap operation, using
the Newton-based method of \cite{AfKlinteberg2018}.

\subsubsection{On-boundary evaluation}

For the on-boundary evaluation outlined in
\cref{sec:boundary-evaluation}, we apply the kernel-split quadrature
for the $\log(\tau-z)$ singularity, using the $n=16$ discretization
points. For a source panel $\Gamma_ i$ this is necessary for target
points $z$ on $\Gamma_i$, and for target points $z$ on neighboring
panels $\Gamma_{i \pm 1}$ such that
$\gamma_i^{-1}(z) \in [-1.4, 1.4]$. Note that we do not need a Newton
solve to find $t^* = \gamma_i^{-1}(z)$, since we from the
discretization know $t$ such that $\gamma_{i \pm 1}(t) = z$ (it is one
of the Gauss-Legendre nodes on $[-1,1]$), and $\gamma_i$ and
$\gamma_{i \pm 1}$ describe the same function, differing only by a
known linear scaling of the argument,
$\gamma_i(t) = \gamma_{i \pm 1}(at+b)$.

\subsubsection{Near-boundary evaluation}
\label{sec:impl-near-bound-eval}

For the near-boundary evaluation of
\cref{sec:near-bound-eval,sec:gradient}, we determine which
point-panel pairs to evaluate using kernel-split quadrature in a
series of steps. For a given target point $z$, we first determine
which panels are within $1.2 h$ of $z$, where $h$ is the panel
length. This is implemented by first sorting the boundary points into
bins of size $h_{max}$ (the maximum panel length on $\bdry$), and then
computing the distances to all boundary points in bins adjacent to
that of the target point. For a panel $\Gamma$ satisfying this first
criterion, we next determine the Bernstein radius $\rho$ of the
preimage of $z$, under that panel's parametrization. If
$\rho \ge 3.5$, then the underlying 16-point Gauss-Legendre quadrature
is accurate at $z$, and no special treatment is necessary. Otherwise,
we proceed by first interpolating all the quantities on $\Gamma$ to
the nodes of a 32-point Gauss-Legendre quadrature, using barycentric
Lagrange interpolation \cite{Berrut2004} (we refer to this as
\emph{upsampling} $\Gamma$). We then evaluate the quadrature from
these 32 points using either the Gauss-Legendre weights, or, if
$\rho < \sqrt{3.5}$, weights computed using the kernel-split scheme.
This near-boundary evaluation can be summarized as:
\begin{itemize}
\item If $\rho \ge 3.5$: Use the 16-point Gauss-Legendre
  quadrature on the underlying panel.
\item If $\sqrt{3.5} \le \rho < 3.5$: Use the 32-point Gauss-Legendre
  quadrature on the upsampled panel.
\item If $\rho < \sqrt{3.5}$: Use kernel-split quadrature on the
  upsampled panel.
\end{itemize}
The upsampling to 32-points panels in the near evaluation is useful
for two reasons: It increases the accuracy of the kernel-split
quadrature, and it increases speed by reducing the number of points
where kernel-split quadrature is necessary.

To make near-boundary evaluation efficient when time stepping, we
precompute and store a set of correction weights for each point-panel
interaction pair that is not accurately evaluated using direct
16-point Gauss-Legendre quadrature. This allows us to first compute
the direct quadrature everywhere in $\Omega$ using the FMM of
\cref{sec:fast-summation}, and then use the correction weights to
subtract off the direct quadrature and add the corrected quadrature
for the points that are close to the boundary. The cost of this
application is very small in relation to the cost of the FMM. The
required storage is also small, since the precomputed weights for
evaluating both layer potential and gradient only requires storage of
$6 \times 2 \times 16$ values for each point-panel pair. In addition,
most near evaluation points interact with only one or two panels, and
we only need to do this for a small band of points close to $\bdry$.

For the above near-boundary evaluation scheme to be robust for any
target point in $\domain$, two special cases must be considered:

\begin{itemize}
\item The first special case occurs when a target point is very close to one of
the discretization nodes on $\bdry$. Then, first adding the direct
contribution through the FMM, before subtracting it off again in the
correction weights, introduces a cancellation error. This error is
$\mathcal O(r^{-1})$ in $\uu$ and $\mathcal O(r^{-2})$ in $\nabla\uu$,
where $r$ is the distance between the target point and the nearest
boundary point. To avoid this, we modify the FMM to ignore direct
interactions when $r < 0.05h$, where $h$ is the length of the panel
to which the source point belongs.

\item The second special case occurs when a target point is very close to an edge
between two panels. In this case, cancellation errors occur in the
computation of the exact integrals $p_k$ \eqref{eq:specquad_pk} used
in the complex interpolatory quadrature. To avoid this, we use a
variation of the panel-merging strategy suggested in
\cite{Helsing2008}. The two neighboring 16-points panels are merged,
using interpolation of boundary data, into a temporary 32-point panel,
from which the quadrature is evaluated as described above. This is
done when the distance $r$ between target point and panel edge
satisfies $r < 0.05 \min(h_1, h_2)$, where $h_1$ and $h_2$ are the
lengths of the two nearby panels.
\end{itemize}

\subsubsection{Dependence on the parameter $\alpha$}
\label{sec:largealpha}

Recall that the modified Stokes equation
\eqref{eqn:modstokes_homogeneous} depends on the parameter
$\alpha=\sqrt{\Re / \deltat}$, which we can expect to vary by orders
of magnitude for different values of $\Re$ and $\deltat$. For our
method to be robust, our quadrature scheme needs to be able to
accurate integrate the stresslet $T$ for a wide range of
$\alpha$. Considering the explicit form \eqref{eq:stresslet_explicit},
having large or small values of $\alpha$ does not appear to be a major
difficulty, as it mainly rescales the action of the kernel. However,
even though the stresslet functions $\mathcal T_1$--$\mathcal T_3$
decay algebraically with $z=\alpha\norm{\xx-\yy}$, the functions
$\mathcal T_1^S$--$\mathcal T_3^L$ into which they are decomposed in
\cref{eq:T1_split,eq:T2_split,eq:T3_split} grow exponentially with $z$
(see \cref{fig:kernel_split}). The root cause of this can be found in
the asymptotic result $I_0(z) \sim I_1(z) \sim e^{z}/\sqrt{2 \pi z}$
when $z \to \infty$ \cite[\S10.30]{NIST:DLMF}. From a numerical point
of view, this is troublesome. Recall that the kernel-split quadrature
of \cref{sec:near-bound-eval}, which we use for the logarithmic
kernel, works by rewriting a convolution integral as
\begin{align}
  \int f(z) K(\tau-z) \dif\tau =
  \int f(z) K^S(\tau-z) \dif\tau + \int f(z) K^L(\tau-z) \log(\tau-z) \dif\tau,
  \label{eq:schematic_log_split}
\end{align}
and then evaluating the two integrals on the right hand side using
Gauss-Legendre quadrature and complex interpolatory quadrature,
respectively. Recall also that both quadratures rely on the remaining
integrand ($f(z) K^S(\tau-z)$ and $f(z) K^L(\tau-z)$, respectively)
being well approximated by a polynomial. Now, when $K^S$ and $K^L$ are
of opposite sign, and both grow as $e^{\alpha|\tau-z|}$, this scheme
will deteriorate quickly if $\alpha|\tau-z|$ gets large. This is due
to a combination of two factors: interpolation errors, because the
exponential factor is no longer well approximated by a polynomial of
order $n-1$, and catastrophic cancellation, due to the large
magnitudes of the terms in the quadrature.  This effectively
introduces a maximum allowable value for the exponent
$\alpha|\tau-z|$. Since kernel-split quadrature is only applied to
target points $z$ that are $\mathcal O(h)$ close to the source panel,
this translates into a maximum allowable length of the source panels
in the kernel-split quadrature. This can be expressed as the following
criterion, which has been empirically determined for our problem,
\begin{align}
  \alpha h \le 4.5 .
  \label{eq:alpha_h_crit}
\end{align}

\begin{figure}
  % Generated by: julia/output_scripts/quad_geo.jl
  \centering
  \includegraphics[height=0.3\textwidth]{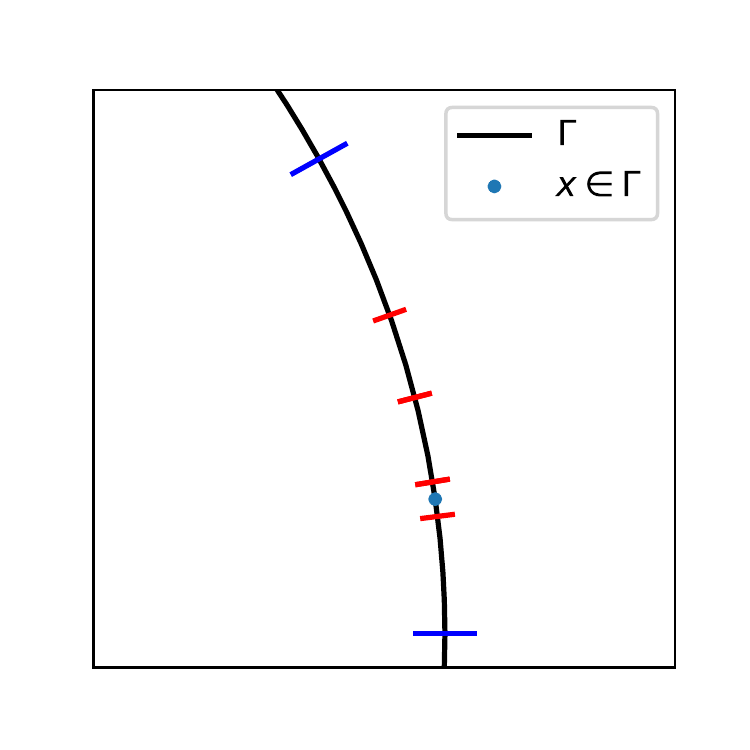}  
  \hspace{2em}
  \includegraphics[height=0.3\textwidth]{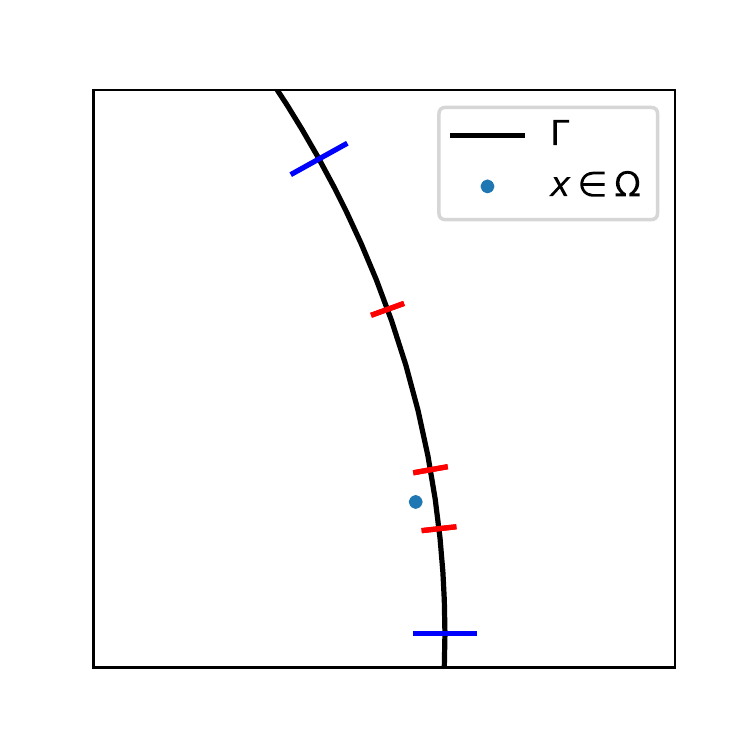}
  \caption{Illustration of the scheme used to subdivide the source
    panel when $\alpha$ is large, for a target point $\xx$ either on
    or near the boundary $\Gamma$. The blue markings indicate the
    boundaries of the underlying source panel $\Gamma$, while the red
    markings indicate the boundaries of the new subpanels introduced
    when evaluating the layer potential at $\xx$.}
  \label{fig:subdivision}
\end{figure}

While we could ensure that all panels satisfy the above criterion at
discretization, this would force us to resolve the boundary more than
necessary, and ultimately impose a time/space discretization restriction of
the form $h^2/\deltat < C$ on our method. To avoid this, we use the
scheme outlined in \cite{AfKlinteberg2019}, which builds
on the observation that the source panel $\Gamma$ and the source
density $\v\mu$ are well resolved by discretization, and that the
difficulties lie in properly resolving the kernel components
$\mathcal T_1^S$--$\mathcal T_3^L$, which are analytically known
functions.
We can therefore subdivide the source panel $\Gamma$ into a set of
temporary subpanels, from which the contributions at $z$ are
evaluated, using boundary data interpolated from $\Gamma$. This allows
us to choose the subpanels such that the ones requiring kernel-split
quadrature, based on their relative distance from $z$, also satisfy
\eqref{eq:alpha_h_crit}. An efficient way of doing this is by using
subpanels that are successively refined in the direction of the target
point $z$, as illustrated in \cref{fig:subdivision}. 
The resulting quadrature scheme is robust, and accurately evaluates
the layer potential both for large $\alpha$ and target points very
close to $\bdry$. The hierarchical structure of the subdivision
algorithm makes the scheme fast, even for large $\alpha$. In addition,
it only incurs an additional cost at the precomputation step; the
action of the quadrature over the subpanels is reduced to a set of
modified quadrature weights at the original 16 sources nodes, using
interpolation matrices.

\subsubsection{Identifying interior points}
\label{sec:ident-inter-points}

An important component of our method is knowing which of the grid
points in $B$ lie in $\domain$, as illustrated in
\cref{fig:embedded_domain}. There are many different ways of
determining this, here we choose a layer potential approach. We first
evaluate the Laplace double layer identity
\begin{align}
  \phi(\xx) = -\frac{1}{2\pi} \int_\bdry \frac{(\xx-\yy)\cdot\v\nhat}{\norm{\xx-\yy}^2} \dSy =
  \begin{cases}
    1, \quad \xx \in \domain,\\
    0, \quad  \xx \notin \domain,
  \end{cases}
  \label{eq:laplace-ident}
\end{align}
using an FMM and our underlying 16-point Gauss-Legendre quadrature. In
the first pass, any point such that $\phi(\xx)>\frac{1}{2}$ is marked
as being interior. In the second pass, we use the algorithm of
\cref{sec:impl-near-bound-eval} to find any points that are too close
to $\bdry$ for the quadrature approximation of
\eqref{eq:laplace-ident} to be accurate. For these points, we find the
closest boundary point $\xx_c$ and mark them as interior if
$(\xx-\xx_c)\cdot\v\nhat(\xx_c) > 0$. This can still fail if
$\norm{\xx-\xx_c}$ is smaller than the distance between $\xx_c$ an its
neighboring points on $\bdry$. In such cases we compute the preimage
of $\xx$ under the parametrization of the nearest panel, as described
in \cref{sec:determ-where-apply}, and mark the point as interior if
$\imag[\gamma^{-1}(x_1+ix_2)] > 0$. This algorithm identifies interior
points in a fast and robust fashion, and builds on methods that are
already present in our code.

%%% Local Variables: 
%%% TeX-master: "nseIE"
%%% End: 

%%%%%%%%%%%%%%%%%%%%%%%%%%%%%%%%%%%%%%%%%%%%%%%%%%%%%%%%%%%%%%%%%%%%%%%%
\section{Numerical examples}
\label{sec:numerical-examples}

\subsection{Implementation}

The software implementation of our scheme is for the most part written
in Julia \cite{Bezanson2017}, and is available for reference as open
source code \cite{inse-fiem-2d}.  The implementation depends on
several external packages: For function extension, we use a Matlab
implementation of PUX written by the authors of \cite{Fryklund2018},
available at \cite{puxdemo}. This in turn makes use of the RBF-QR
algorithm \cite{Fornberg2011}, as implemented in
\cite{rbfqr}. Non-uniform FFTs are computed using the library FINUFFT
\cite{Barnett2018}. Our FMM implementation is written in Fortran, and
is based on FMMLIB2D \cite{fmmlib2d}. We also use FMMLIB2D directly
for identifying interior points (see
\cref{sec:ident-inter-points}). The fast direct solver is written in
Julia, but is based on a Matlab implementation of the Stokes solver
used in \cite{Marple2016}, provided by A. Gillman. The low rank matrix
approximations used in the fast direct solver are computed using the
interpolative decomposition (ID) \cite{Liberty2007,cheng_2005}, as
implemented in LowRankApprox.jl \cite{Ho2018}. Our Julia code is
single-threaded, though several of the libraries used, including BLAS
operations, are multi-threaded.

The scheme is designed to be as efficient as possible in solving the
inhomogeneous problem \eqref{eq:modstokes} for a fixed $\alpha$, which
is equivalent to taking IMEX or SISDC steps with a fixed time step
size. To achieve this, we precompute the fast direct solver, the
quadrature weights required to evaluate the layer potential at the
near-boundary grid points, the FMM tree, and the main system matrix of
PUX. This means that there is a significant cost associated with
updating either the discretization or time step size. We therefore run
all of our simulations with a constant time step size, and let SISDC
have equisized substeps.

Throughout the below tests, we denote by $\uu$ the reference solution,
and by $\tilde\uu$ the computed solution. For a discretization with
grid points $\{\xx_j\}_1^\Ndom$, we measure the error $\v e = \tilde\uu-\uu$ in the $\ell^2$
and $\ell^\infty$ norms, computed as
\begin{align}
  \norm{\v e}_2
  = \sqrt{\sum_{j=1}^\Ndom \left( e_1(\xx_j)^2 + e_2(\xx_j)^2 \right)}
  \quad \mbox{ and } \quad
  \norm{\v e}_\infty
  =  \max_{\substack{1 \le j \le \Ndom \\ i=1,2}} \abs{e_i(\xx_j)}. 
\end{align}

Where reported, the computation times are from a computer with 32 GB
of memory and an Intel Core i7-7700 CPU (4 physical cores, 3.6 GHz base
frequency).

\subsection{Homogeneous stationary problem}
\label{sec:homog-stat-probl}

\begin{figure}
  % Generated by: julia/output_scripts/convergence_homogeneous.jl
  \centering
  \includegraphics[width=.4\textwidth]{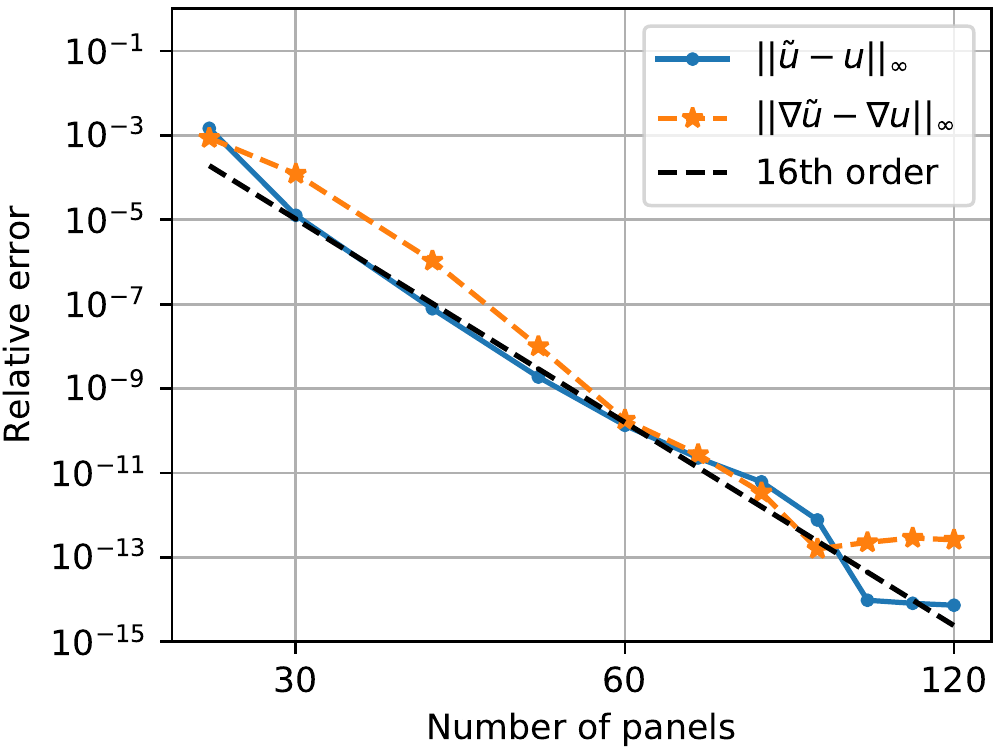}
  \hspace{.1\textwidth}
  \includegraphics[width=.4\textwidth]{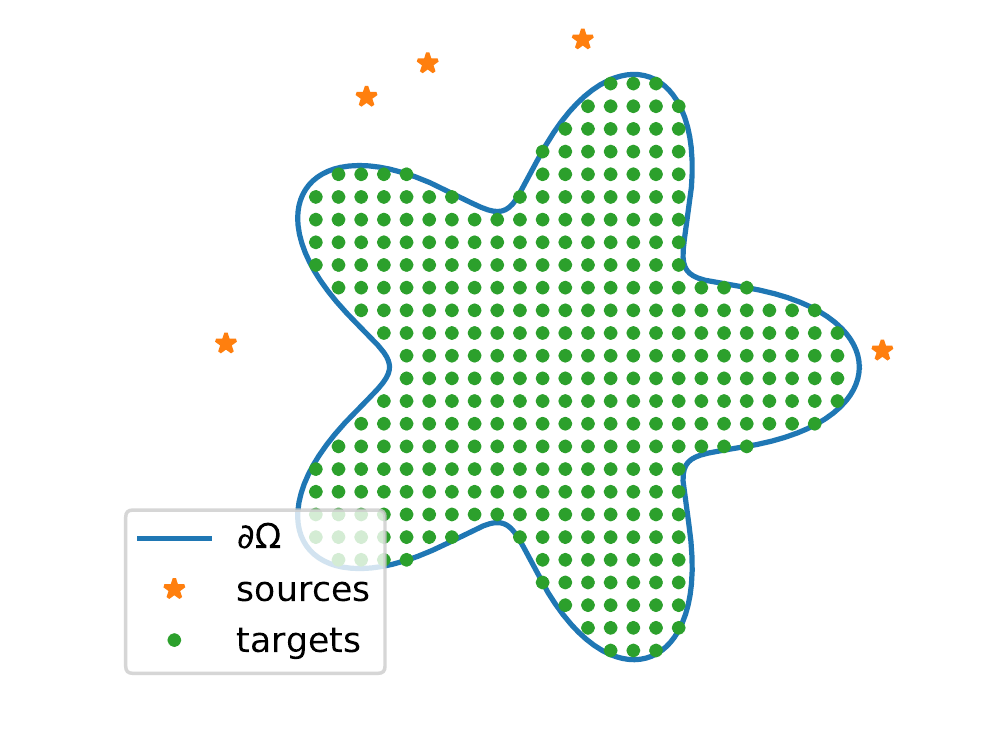}
  \caption{Convergence for the homogeneous solver, with
    $\alpha=10$. The solution is given by stokeslets at the source
    points, and we measure the maximum relative errors in the solution
    and the gradient at the target points.}
  \label{fig:convergence_homogeneous}
\end{figure}

To validate the solver for the homogeneous equation
\eqref{eqn:modstokes_homogeneous}, described in
\ref{sec:comp-homog-corr}, we set up the following problem: The
boundary $\bdry$ is the starfish described by the parametrization
\begin{align}
  \v g(t) = R \left( 1 + A \cos 5 t \right)
  \left(
  \begin{array}{c}
    \cos t \\ \sin t
  \end{array}
  \right), \quad t \in [0, 2\pi) ,
  \label{eq:starfish}
\end{align}
with $R=1$ and $A=0.3$. The Dirichlet boundary condition is given by a
sum of 5 stokeslets,
\begin{align}
  \uu(\xx) = \ff(\xx) = \sum_{n=1}^5 S(\xx, \yy_n) \v q_n, \quad \xx \in \bdry.
\end{align}
This has the exact solution $\uu=\ff$ in $\domain$.  The source
locations $\yy_n$ are randomly placed on a circle of radius 1.4, and
the strengths $\v q_n$ are drawn randomly from $\mathcal U(-1,1)$. We
set $\alpha=10$, and discretize $\bdry$ using a varying number of
panels, from 30 to 130. For each discretization, we first solve the
integral equation by forming the dense linear system and solving it
directly, and then evaluate the solution at $358$ uniform target points
in $\domain$ using our fast summation and near-boundary
quadrature. The geometry of the problem and the convergence results
are shown in \cref{fig:convergence_homogeneous}. We expect to see a 16th
order convergence, since that is the order of the Gauss-Legendre
panels, and that is also what we observe, both in $\uu$ and
$\nabla\uu$.

\subsection{Inhomogeneous problem convergence}
\label{sec:inhom-probl-conv}

\begin{figure}
  \centering
  % Plots below have been trimmed with `mogrify -trim <name>.png`
  % Generated by: julia/output_scripts/convergence_full.jl    
  \includegraphics[width=.3\textwidth]{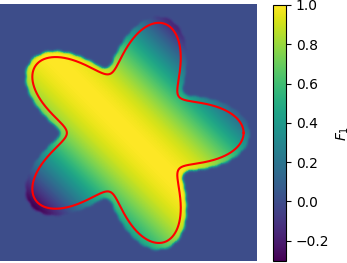}
  \hspace{0.03\textwidth}
  % Generated by: julia/output_scripts/alpha_convergence.jl      
  \includegraphics[width=.3\textwidth]{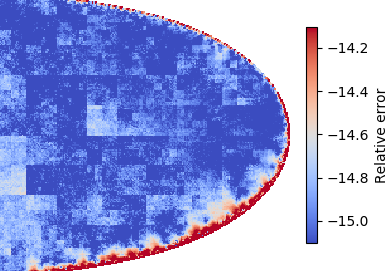}
  \hspace{0.03\textwidth}
  \includegraphics[width=.3\textwidth]{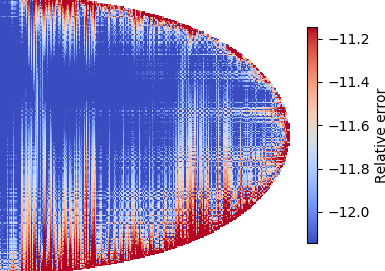}
  \caption{Left: Function extension of the first component of
    \eqref{eq:full_rhs}, using PUX with partition radius $R=0.15$.
    Center and right: Zoom of the rightmost starfish arm, showing the
    error in $\uu$ (pointwise maximum) for $N=800$ and $\alpha=1$
    (center) and $\alpha=1000$ (right). Note that the color scale is
    different between the two pictures. At $\alpha=1$, the box-like
    structure in the error is due to the FMM tree. The error near the
    boundary is due to PUX; the scale is too large to be due to the
    layer potential near evaluation.  The structure in the error for
    $\alpha=1000$ is typical for errors in high Fourier modes, an
    effect of the reduced smoothing for large $\alpha$.}
  \label{fig:convergence_full_fields}
\end{figure}

Next, we validate our composite solver for the full stationary problem
--- the inhomogeneous modified Stokes equations
\eqref{eq:modstokes}. We again use the starfish geometry
\eqref{eq:starfish}, and let the right hand side $\FF$ be a simple
oscillation,
\begin{align}
  \FF(\xx) = (-1, 2) \cos (x_1+x_2), \quad \xx \in \domain .
  \label{eq:full_rhs}
\end{align}
On a periodic box with sides $2\pi$, this has an exact solution given
by a single Fourier mode \eqref{eq:uP_Fourier},
\begin{align}
  \uu_{\text{per}}(\xx) = (3, -3) \frac{\cos(x_1+x_2)}{2(2 + \alpha^2)}, \quad \xx \in \domain .
\end{align}
By setting the Dirichlet boundary condition $\uu=\uu_{\text{per}}$ on
$\bdry$, we get a problem with exact solution $\uu=\uu_{\text{per}}$
in $\domain$. This tests the complete solver setup, since $\FF$ is
only given inside $\domain$, and then smoothly extended into the
bounding box $B$ using PUX (see
\cref{fig:convergence_full_fields}). We discretize $\bdry$ using 400
panels, and set the PUX radius to $R=0.15$. In $\domain$ we use a uniform
grid with $N \times N$ points in the smallest square box that bounds
$\bdry$.

\begin{figure}
  \centering
  % Generated by: julia/output_scripts/convergence_full.jl  
  \includegraphics[width=.4\textwidth]{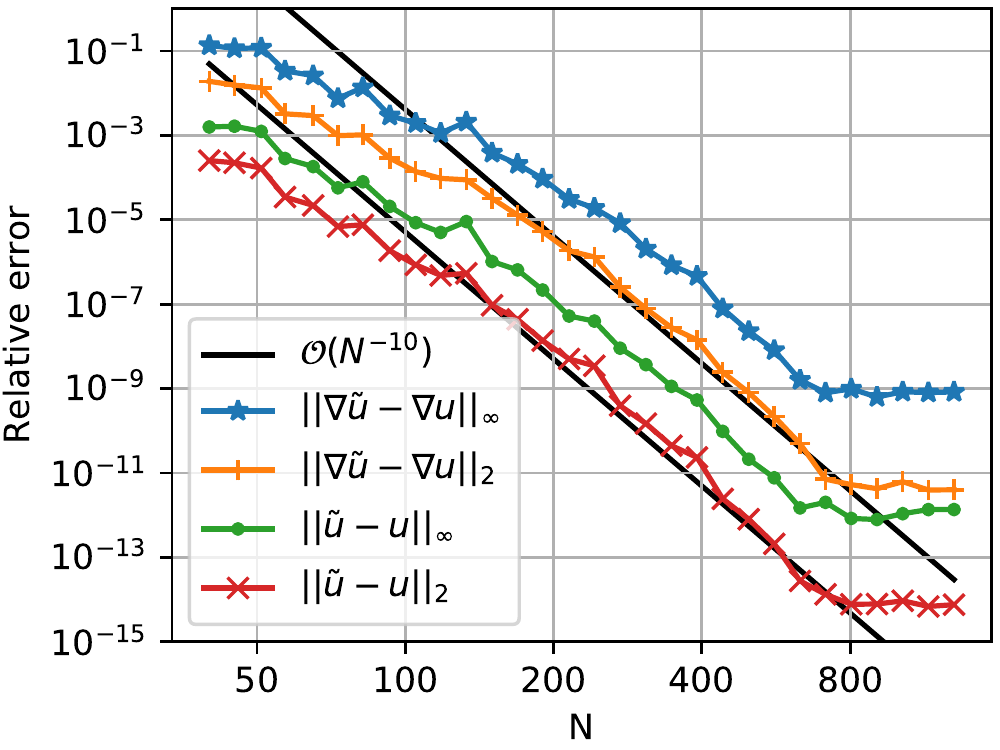}
  \hspace{.1\textwidth}
  % Generated by: julia/output_scripts/alpha_variation.jl  
  \includegraphics[width=.4\textwidth]{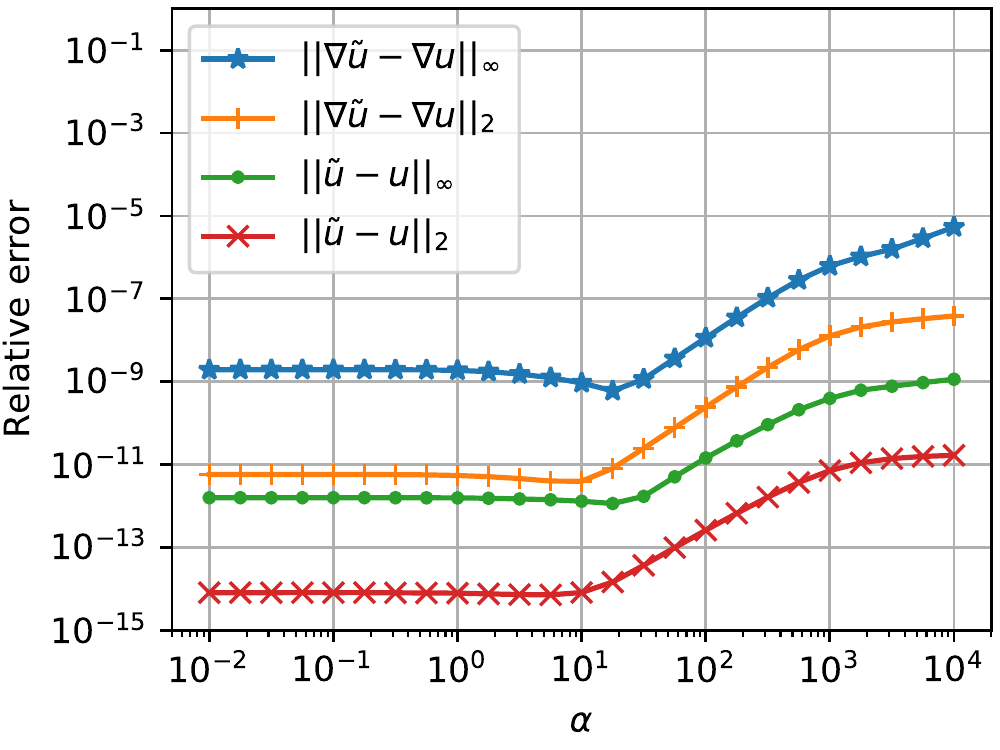}  
  \caption{Left: Convergence for full problem, $\alpha=10$. Right: Varying $\alpha$ when $N=800$.}
  \label{fig:convergence_full}
\end{figure}

To test convergence in the volume grid, we set $\alpha=10$ and solve
the above problem for a wide range of $N$, from 40 to 1300. The
results, shown in \cref{fig:convergence_full}, indicate that we have
10th order convergence in both $\uu$ and $\nabla \uu$. This is in
accordance with the results in \cite{Fryklund2018}. Note that the
errors differ by about two orders of magnitude between 2-norm and
max-norm, and also by about two orders of magnitude between $\uu$ and
$\nabla\uu$. A maximum relative error in $\nabla\uu$ around $10^{-9}$
appears to be the best that we can achieve for this particular
problem. This is due to errors in the solution close to the boundary,
see \cref{fig:convergence_full_fields}.

To test the dependence on the parameter $\alpha$, defined as
$\alpha^2={\Re/\deltat}$, we also solve the above problem with
$\alpha$ varying by several orders of magnitude. This time we keep $N$
fixed at $N=800$, where the solution was fully converged for
$\alpha=10$. The results, also shown in \cref{fig:convergence_full},
clearly show that accuracy suffers when $\alpha$ gets very large.  We
believe that this is error originates in the particular solution, as
the homogeneous solver can handle a very wide range of $\alpha$, due
to the scheme outlined in \cref{sec:largealpha}.~\ref{sec:ident-inter-points} Specifically, the
smoothing effect of the Fourier multiplier
\eqref{eq:periodic_stokeslet} is reduced for large $\alpha$, which
exposes any irregularities in the function extensions. This problem is
even more pronounced in the computation of the gradient
\eqref{eq:uP_Fourier_deriv}, which helps explain the larger errors in
$\nabla\uu$.

\subsection{Time convergence}

To validate our complete solver for the Navier-Stokes equations
\eqref{eq:navier-stokes}, we consider the problem of viscous spin-down
in a cylinder: A cylinder of radius $a$ is filled with viscous fluid,
and at $t=0$ both the cylinder and the fluid are rotating with angular
velocity $\Omega$. Then (for $t>0$) the cylinder suddenly stops
rotating, while the fluid keeps rotating until it comes to rest. In
polar coordinates,
\begin{align}
  \begin{aligned}
    \uu(r,t) &= u_\theta(r,t) \v{\hat\theta}, && \\
    u_\theta(r,0) &= \Omega r, & r &\le a, \\
    u_\theta(a, t) &= 0, & t &> 0.
  \end{aligned}
\end{align}
The exact solution to this problem is known \cite[p.45]{Acheson1990},
\begin{align}
  u_\theta(r, t) = -2 \Omega a \sum_{n=1}^\infty
  \frac{J_1(\lambda_n r / a)}{\lambda_n J_0(\lambda_n)}
  e^{-\lambda_n^2 t / (a \Re)} .
\end{align}
Here $J_k$ is the Bessel function of order $k$, and $\lambda_n$ is the
$n$th positive root of $J_1$.

\begin{figure}
  % Generated by: julia/dev/demo_spindown.jl
  \centering
  \includegraphics[width=0.9\textwidth]{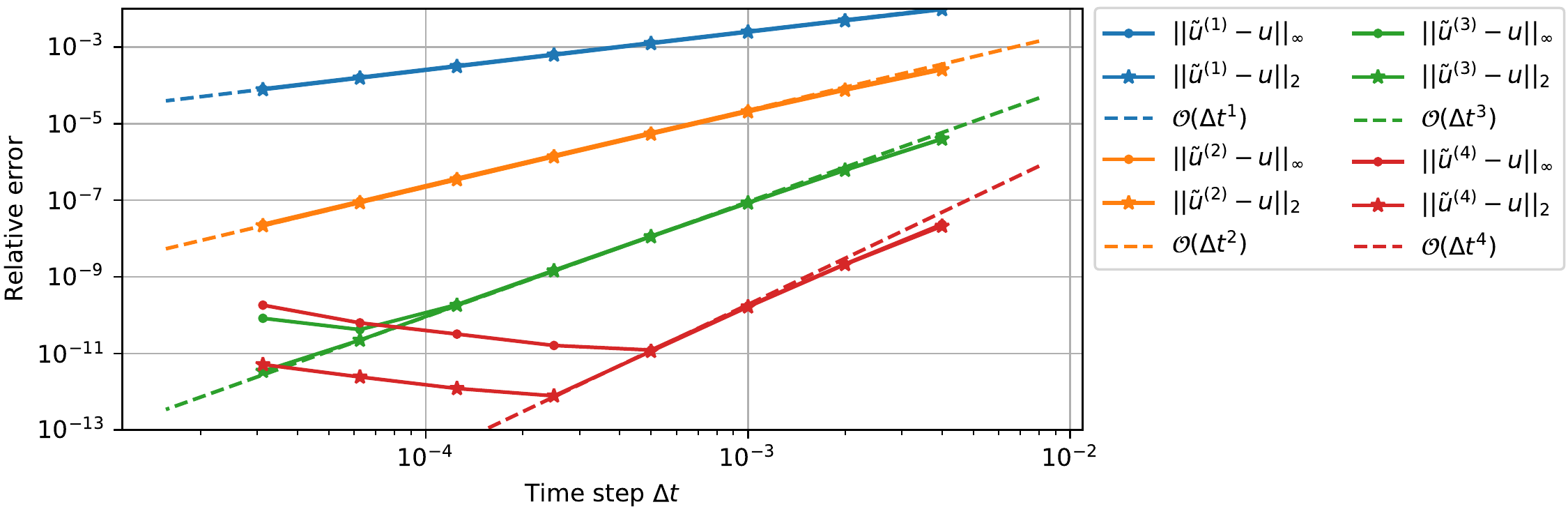}
  \caption{Time convergence for spin down problem, $\tilde u^{(K)}$
    denotes solution computed using SISDC order $K$. For small
    $\Delta t$, the lowest attainable error can be seen to increase
    with decreasing $\Delta t$. This is related to the results in
    \cref{fig:convergence_full}, which show that the error grows with
    increasing $\alpha$.}
  \label{fig:convergence_time}
\end{figure}

We solve the above problem for $a=\Omega=\Re=1$. In order to avoid the
time-discontinuity in the boundary condition at $t=0$, we initialize
our solver at $t=0.05$. The flow field is then sufficiently smooth, so
that we can observe high-order convergence. We run our simulation on
the interval $t \in [0.05, 0.054]$, with an initial time step length
$\Delta t=0.002$, which is then successively halved.  Our
discretization uses a $500 \times 500$ volume grid, 200 panels on the
boundary, and PUX radius $R=0.3$. We time step using SISDC of order 1
through 4, and compare the results in
\cref{fig:convergence_time}. Convergence is as expected in both
$\ell^2$ and $\ell^\infty$, although the fourth order convergence only
just reaches the expected rate before being overtaken by what appears
to be an error that grows with decreasing $\Delta t$, and is
approximately the same for $K=3$ and $K=4$. To explain this, recall
that with time step length $\Delta t$ and SISDC order $K$, the IMEX
substep length is $\deltat = \Delta t / (K-1)$, such that
\begin{align}
  \alpha = \sqrt{\frac{(K-1) \Re}{\Delta t}} .
\end{align}
For this time convergence test, the smallest time steps used
correspond to values $\alpha$ in the region where the errors in the
stationary solver start increasing, as seen in
\cref{fig:convergence_full}. This explains the increasing error, and
also why it is approximately the same between $K=3$ and $K=4$.

\subsection{Stability}
\label{sec:stability}

\begin{figure}
  % Generated by:
  % julia/dev/time_stability.jl
  \centering
  \subfloat[][\\$\deltat=0.04$\\$\Re=50\quad\:$]{  
    \includegraphics[height=0.2\textwidth]{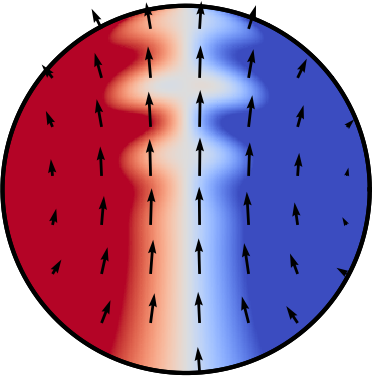}
    \label{fig:stability_1}
  }
  \subfloat[][\\$\deltat=0.02$\\$\Re=50\quad\:$]{    
    \hspace{0.2cm}
    \includegraphics[height=0.2\textwidth]{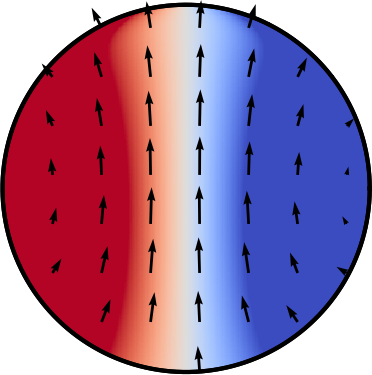}
    \label{fig:stability_2}
  }
  \subfloat[][\\$\deltat=0.02$\\$\Re=100\;\;$]{        
    \hspace{0.2cm}  
    \includegraphics[height=0.2\textwidth]{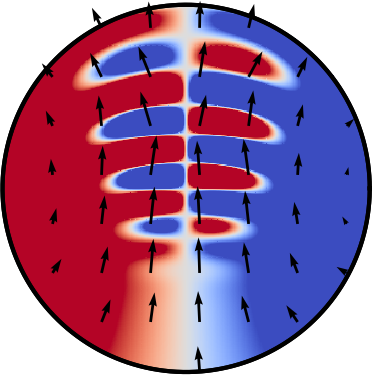}
    \label{fig:stability_3}        
  }
  \subfloat[][\\$\deltat=0.01$\\$\Re=100\;\;$]{        
    \hspace{0.2cm}  
    \includegraphics[height=0.2\textwidth]{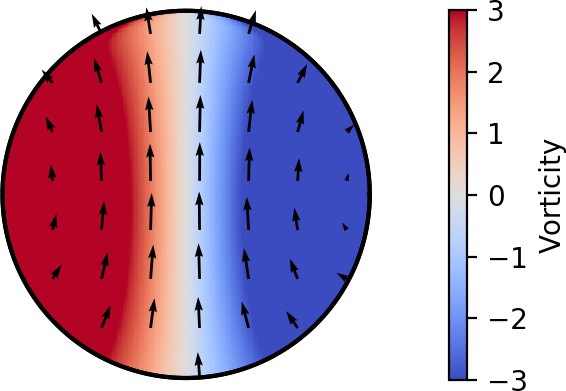}
    \label{fig:stability_4}    
  }
  \caption{
    Simulations of the flow \eqref{eq:instab_flow}, with varying $\deltat$ and $\Re$.
    %.
    Cases \protect\subref{fig:stability_2} and
    \protect\subref{fig:stability_4} satisfy the heuristic stability
    condition \eqref{eq:stab_cond}, and develop a steady state. Snapshots
    shown are after 200 time steps.
    % .
    On the other hand, cases \protect\subref{fig:stability_1} and
    \protect\subref{fig:stability_3} violate the stability condition,
    and develop an instability pattern that is more pronounced for
    larger $\Re$. The snapshots shown are after 45 time
    steps; the solution eventually blows up completely if run longer.
  }
  \label{fig:stability}
\end{figure}

To investigate the stability properties of our method, we consider the
simple case of flow in a circle of diameter 1, with boundary condition
\begin{align}
  \uu = -%V
  \sin(\phi)\v\nhat, \quad \xx \in \bdry,
  \label{eq:instab_flow}
\end{align}
where $\phi=\tan^{-1}(x_2/x_1)$ is the angle of the point on the
boundary. Simulating this flow with the IMEX scheme, we observe
instabilities in the form of waves in the vorticity, originating at
the inflow, as shown in \cref{fig:stability}. Through numerical
experimentation, we find that these instabilities occur unless we
satisfy a condition of the form
\begin{align}
  %V^2
  \deltat \Re \le C ,
  \label{eq:stab_cond}
\end{align}
where, for this particular flow problem, $C \approx 1$. The
instabilities do not appear to depend on the spatial scale (radius of
the circle), nor the spatial discretization. The examples in
\cref{fig:stability} were computed using a high-resolution spatial
discretization, but refining or coarsening that discretization has no
effect, as long as it is sufficiently fine for accurately solving the
modified Stokes equation. Switching to higher-order SISDC is
stabilizing, effectively increasing $C$, but does not change the form
of the stability condition.

As an explanation model for our observations, we consider the standard
stability analysis of the 1D advection-diffusion equation, which we write as
\begin{align}
  u_t + %V
  u_x - \frac{1}{\Re} u_{xx} = 0.
\end{align}
This serves as a linearized model of the non-dimensionalized
Navier-Stokes equations.
% , for a flow with maximum velocity $V$.
We discretize this using first-order
IMEX in time and a Fourier series in space, writing
$u(x, t_n) = \sum_k \hat u_k^n e^{ikx}$. This is representative of how
we solve the particular problem \eqref{eqn:modstokes_particular}, and
diagonalizes our model problem to the difference equation
\begin{align}
  \frac{\hat u_k^{n+1} - \hat u_k^{n}}{\deltat} + ik%V
  \hat u_k^{n} + \frac{1}{\Re} k^2 \hat u_k^{n+1} = 0 .
\end{align}
For this to have a stable solution, it must hold that
$|\hat u^{n+1}_k| \le |\hat u^{n}_k|$ for all $k$. It is
straightforward to show that this leads to the condition
\begin{align}
  % V^2
  \deltat \Re \le 2 .
\end{align}
This is consistent with our observed condition
\eqref{eq:stab_cond}. Although derived using a simplified model, it is
a strong indication of the origin of the observed instability.

Considering the physics of our problem, our results are perhaps not
surprising. The smaller the viscosity is compared to the velocity, the
more advection dominated the flow is. It is then only natural that the
timestep must be small for the solution to be stable when the advective
term is treated explicitly.

\subsection{Flow past obstacles}
\label{sec:flow-past-obstacles}

\begin{figure}
  \centering
  % Generated by: julia/dev/obstacles_plot.jl
  \includegraphics[width=0.75\textwidth]{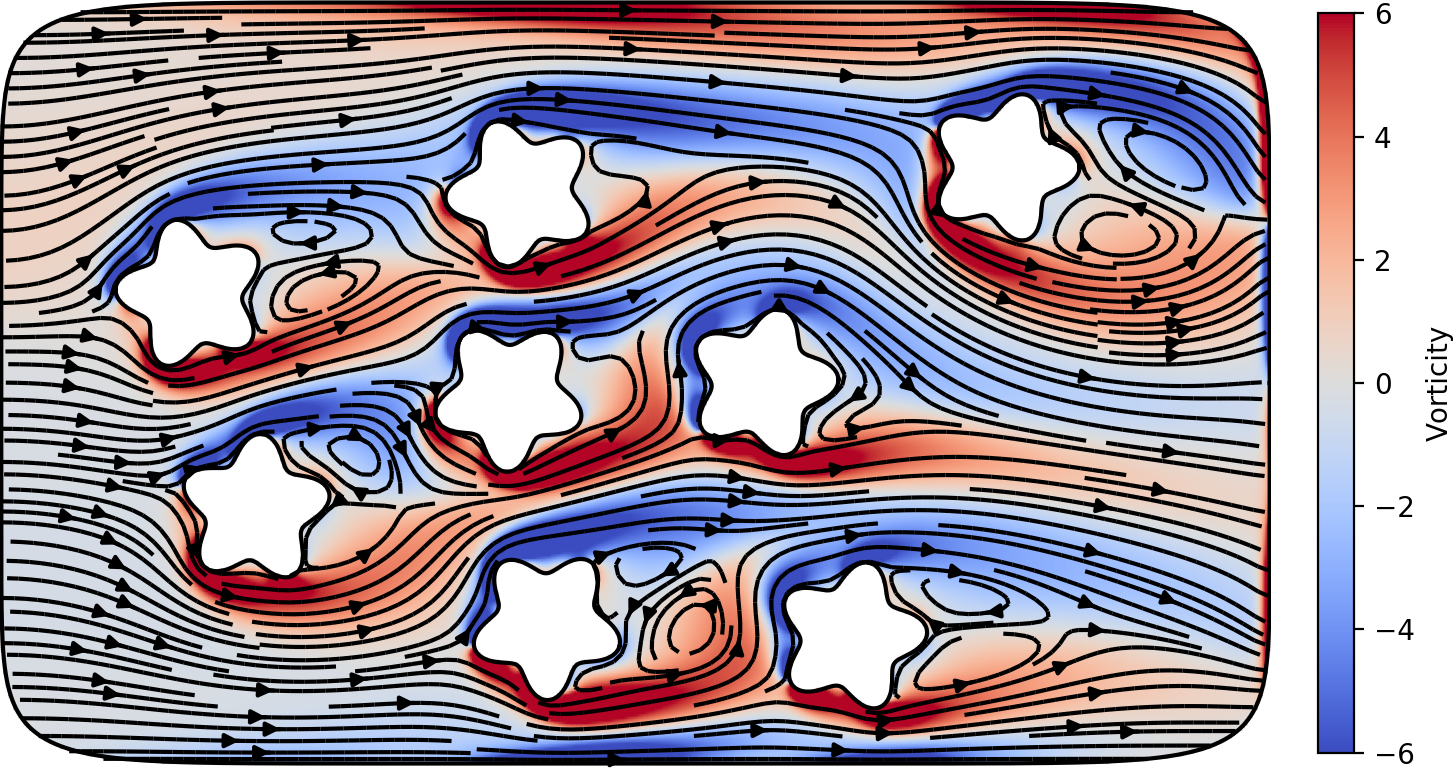}
  \caption{Streamlines and vorticity field of flow past a collection
    of starfish at $\Re=30$.}
  \label{fig:demo}
\end{figure}

As a demonstration of flow through a relatively complex geometry, we
set up the geometry shown in \cref{fig:demo}. The outer boundary is the
box-like domain described by the parametrization
\begin{align}
  \v g(t) = \frac{\left(W \cos t, H \sin t\right)}{\left((\cos t)^p + (\sin t)^p \right)^{1/p}},
  \quad t \in [0, 2\pi),
\end{align}
with $(W,H)=(5,3)$ and $p=10$. The eight starfish-shaped inclusions
are rotated and translated instances of the curve given by
\eqref{eq:starfish}, with $R=0.5$ and $A=0.2$. The boundary conditions
are slip on the outer boundary, $\uu=(1,0)$, and no-slip on the
inclusions, $\uu=(0,0)$. The Reynolds number is set to $\Re=30$.

The problem is discretized using a grid of $500 \times 300$ points in
the domain, $500$ panels on the outer boundary, $50$ panels on each
inclusion, PUX radius 0.4, and second order SISDC with $\Delta
t=0.01$. This gives $\alpha=54.8$. The homogeneous problem is solved
using the fast direct solver (FDS), computed using a tolerance of
$10^{-8}$. This spatial discretization gives an error of
$\norm{\uu - \tilde\uu}_\infty = 5 \cdot 10^{-6}$ on the test problem
used in \cref{sec:inhom-probl-conv}.

The precomputation time for this problem is dominated by the FDS,
which takes 160 s to compute. Remaining precomputations that take more
than one second are 6 s for PUX, and 10 s for the nearly singular
quadrature. Once everything is precomputed, each solve takes around
2.5 s, where the average times for the involved algorithms are PUX:
0.5, FFT+NUFFT: 0.15, FDS: 0.03, FMM: 1.7. With second order SISDC (2
solves per step), each time step takes around 5 s.

We run the simulation for 10000 steps, after which the flow has
reached a steady state. In \cref{fig:demo}, we plot the streamlines
and the vorticity, which is defined as
\begin{align}
  \omega = \dpd{u_2}{x} - \dpd{u_1}{y} .
\end{align}

\subsection{Vortex shedding}
\label{sec:vortex-shedding}

\begin{figure}
  % Generated by:
  % julia/dev/demo_unsteady.jl
  % julia/output_scripts/vortex.jl  
  \centering
  \subfloat[][$\Re=25$]{
    \includegraphics[height=0.27\textwidth]{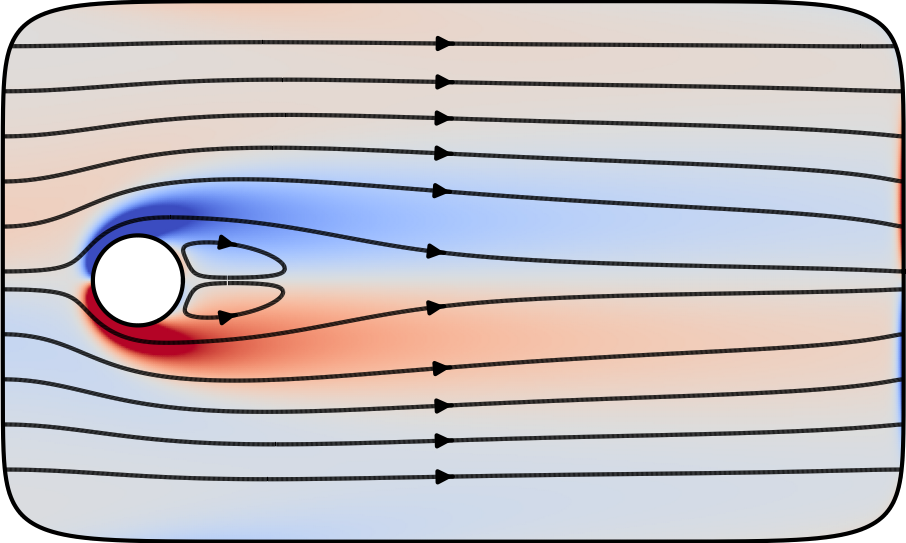}
  }
  \subfloat[][$\Re=50$]{  
    \includegraphics[height=0.27\textwidth]{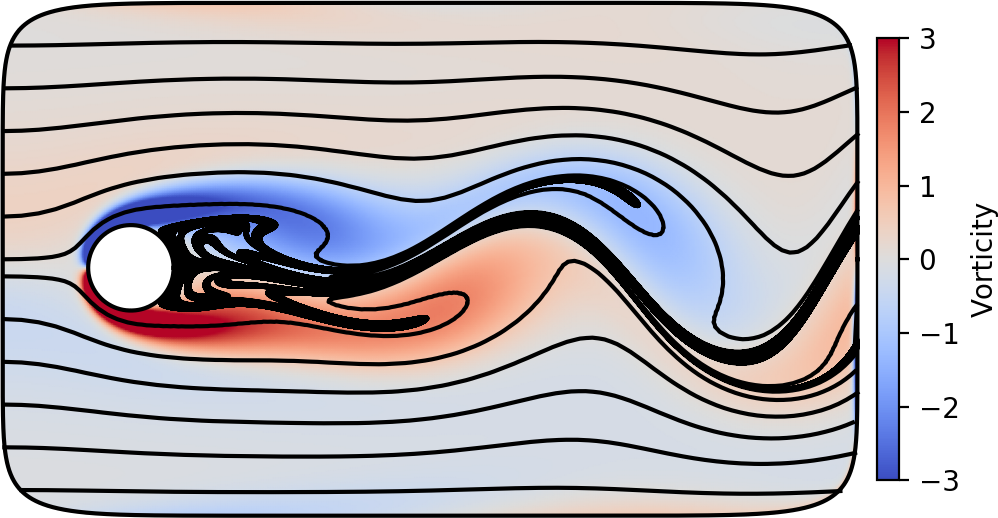}
  }
  \caption{Flow past cylinder, with streaklines drawn in black and
    field colored by vorticity. At $\Re=25$ the flow is steady with
    standing eddies in the wake, while at $\Re=50$ the flow is
    unsteady with a vortex street in the wake.}
  \label{fig:cylinder1}
\end{figure}

\begin{figure}
  % Generated by:
  % julia/dev/demo_unsteady.jl
  % julia/output_scripts/vortex.jl
  \centering
  \subfloat[][$\Re=100$]{
    \includegraphics[height=0.27\textwidth]{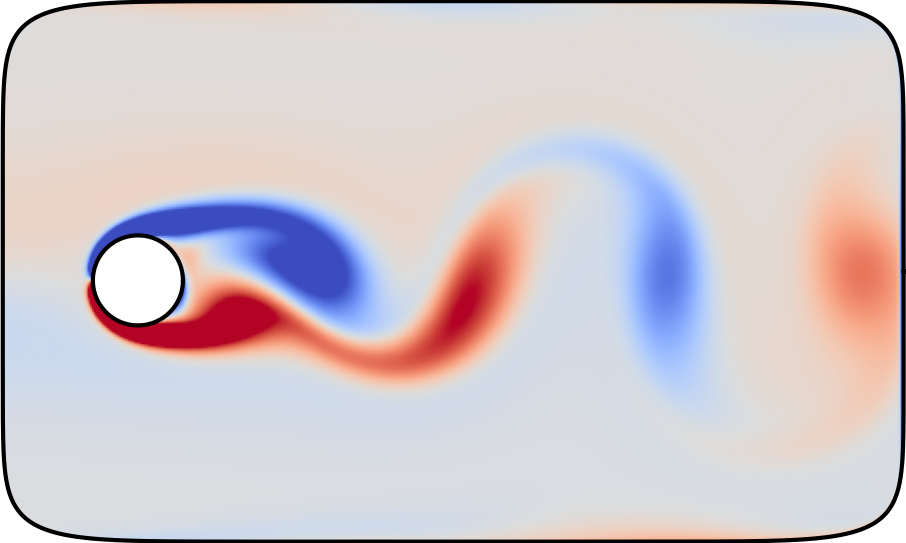}
  }
  \subfloat[][$\Re=200$]{  
    \includegraphics[height=0.27\textwidth]{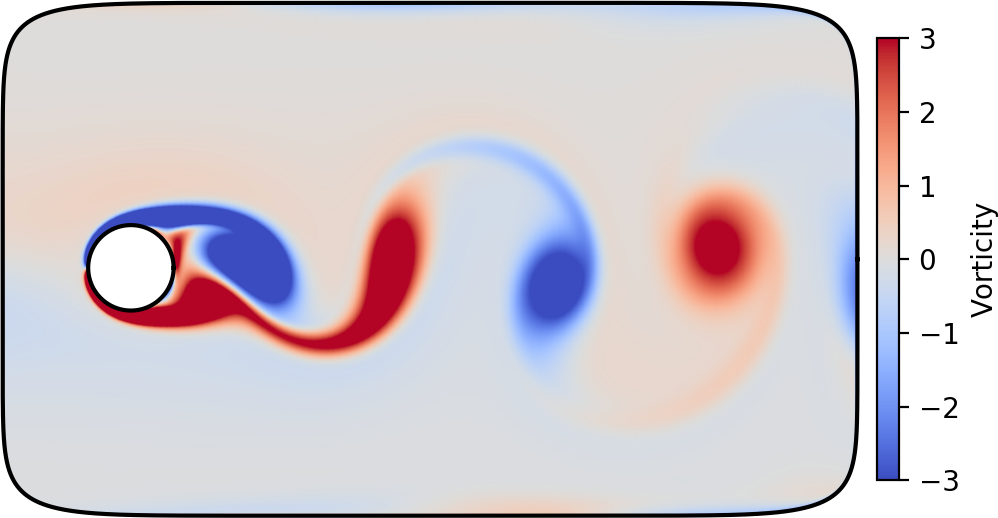}
  }
  \caption{Unsteady flow past cylinder, with a vortex street in the
    wake. Field colored by vorticity.}
  \label{fig:cylinder2}
\end{figure}

To demonstrate the effect of the Reynolds number, we consider the case
of flow past a cylinder for $\Re=\{25,50,100,200\}$, shown in
\cref{fig:cylinder1,fig:cylinder2}. We let the outer domain be the
same box-shaped domain as in \cref{sec:flow-past-obstacles}, with the
same slip boundary condition, and add an inclusion in the form of a
circle with diameter 1 and a no-slip boundary condition. The circle is
centered at $(-3.5, -0.1)$; the vertical offset is added because is
reduces the number of time steps required before the vortex separation
occurs.

We discretize the problem using 500 panels on the outer boundary, and
100 panels on the circle. The domain is discretized using
$500 \times 300$ points and PUX radius $0.5$ for all values of $\Re$
except $\Re=200$, which is discretized using $1000 \times 600$ points
and PUX radius $0.25$, in order to get high accuracy for the larger
$\alpha$. We time step using second order SISDC, and observe that we
need to set $\Delta t= 1/\Re$ for the solution to be stable. This is
consistent with the results in \cref{sec:stability}, and as a result
$\alpha=\Re$. Each solution is timestepped until it is either steady,
or the dynamics of the unsteady solution are fully developed.

The flow behind a cylinder in a free stream develops an oscillating
tail of trailing vortices (a vortex street) for Reynolds numbers above
the critical value, which is around 40. This is consistent with our
simulation, which is steady for $Re=25$, and develops a vortex street
for $\Re=50$ and higher. Even though it is confined to a relatively
small domain, with free-stream boundary conditions at the outer
boundary, the results in \cref{fig:cylinder1,fig:cylinder2} match
qualitatively with other results in the literature,
e.g. \cite{Tafuni2018}.

%%% Local Variables: 
%%% TeX-master: "nseIE"
%%% End: 

%%%%%%%%%%%%%%%%%%%%%%%%%%%%%%%%%%%%%%%%%%%%%%%%%%

\section{Conclusion}
\label{sec:conclusion}

In the preceding we have demonstrated a fast integral
equation method (FIEM) for the Navier--Stokes
equations in general smooth geometries in two dimensions. 
The solver has a number of favorable features.
In terms of accuracy, the solver is 10th-order accurate
in the spatial grid spacing, and we have demonstrated
up to 4th-order convergence in time.
By using an embedded boundary approach combined with a uniform grid,
the method can effortlessly deal with complex geometries. In addition,
the use of potential methods means that the solution satisfies the
incompressibility contraint by construction. This eliminates the need
for projection methods, and the artificial boundary conditions
associated with them. It also makes it straightforward to use
off-the-shelf methods to augment the temporal accuracy of the
underlying IMEX discretization, and one could easily replace the SISDC
scheme with, for example, a linear multistep method based on a
backward differentiation formula (BDF) discretization.

Our numerical results indicate that our method can be used to solve
the unsteady Navier-Stokes equations to high accuracy, for Reynolds
numbers up to at least $\Re=200$. They do however also point towards a
fundamental limitation in using a discretization based on Rothe's
method. For large values of $\alpha$ (corresponding to large Reynolds
numbers or short time steps), we see that more and more spatial
resolution is needed to resolve the underlying kernels; see for
instance, \cref{fig:convergence_full}. This is due to the modified
Stokes equation being fundamentally harder to solve for large
$\alpha$. Combined with the stability restriction
\eqref{eq:stab_cond}, which implies the parameter scaling
$\alpha \sim \Re$, this leads us to believe that our method will
mainly be found useful for flows with low or moderate Reynolds
numbers, in the range of hundreds up to a thousand. 

For computations in two dimensions, there are several straightforward
directions in which the present method can be developed. The integral
equation formulation can be extended to allow for a wider class of
boundary conditions, and to allow post-solution computation of various
quantities, such as the pressure.
While our method applies only to
smooth domains, there has been recent progress in the efficient
discretization of boundary integral equations on domains with corners,
see, for instance, \cite{serkh2016solution, rachh2017solution,
  helsing2008corner, helsing2018integral}. Implementing one of these
schemes would allow accurate computation of flows around geometries
with sharp corners.
It would also be of interest to extend our method to deal with moving
geometries, such as drops, bubbles, vesicles, and particles. Here, the
primary challenge would be to achieve high order accuracy while using
function extension for computing the volume potential.

Although there are several ways in which the present work may be
improved, we believe that the most interesting direction of
development is to extend it to three dimensions (3D). Many of the
tools used in this work extend directly to 3D, but not all. In
particular, the panel-based, kernel-split quadrature is limited to 2D,
but could be replaced by the quadrature-by-expansion (QBX) method
\cite{Klockner2013,Wala2019}. In addition, high performance computing
aspects become much more important in 3D, due to the large amount of
degrees of freedom, and issues of adaptivity and parallelization must
be considered from the start.

%%% Local Variables: 
%%% TeX-master: "nseIE"
%%% End: 

\section*{Acknowledgments}

The authors gratefully acknowledge support from the Knut and
Alice Wallenberg Foundation under grant no. 2016.0410 (LaK), from the
Air Force Office of Scientific Research under grant
FA9550-17-1-0329 (TA), and from the Natural Science and Engineering
Research Council of Canada (MCK).  We also wish to thank Adrianna
Gillman for providing us with an implementation of a fast direct
solver, and Fredrik Fryklund for providing us with an implementation
of PUX.

\clearpage
\appendix
\begin{appendices}
  
\section{Modified Stokes layer potentials}
\label{sec:layer-potentials}

In this section, we provide some technical
results which describe the properties of the
modified Stokes double layer potential and
which explain the use of the nullspace
correction. We closely follow the presentation of
\cite{biros2002embedded}. Details and proofs
may be found in \cite{KimSangtae1991M:pa} for the
Stokes case and \cite{Pozrikidis1992,biros2002embedded}
for the modified Stokes case.

For reference, the homogeneous modified
Stokes equation with the Dirichlet boundary condition
is given by

\begin{align}
  (\alpha^{2} - \Delta)\vv + \nabla p &= 0, 
  \quad &&\xx \in \domain, \label{eq:app_momentum} \\
  \nabla \cdot \vv &= 0, &&\xx \in \domain,
  \label{eq:app_cty}\\
  \vv &= \boldg, &&\xx \in \bdry \; \label{eq:app_dir}.
\end{align}
As in the main text, the divergence-free condition
implies that the boundary data should satisfy

\begin{equation}
  \label{eq:app_compatg}
  \int_\domain \boldg \cdot \v\nhat \, \dSy \; .
\end{equation}

The following energy-type result is convenient
for establishing the uniqueness of various
boundary value problems for the modified Stokes
equation.
\begin{lemma}[Energy \cite{biros2002embedded}]
  \label{lemma:energy}
  Suppose that $\vv$ is satisfies
  \cref{eq:app_momentum,eq:app_cty} on
  a (possibly multiply-connected)
  domain $\domain$ and let $\ssigma(\vv)$ denote
  the associated stress tensor. Then
  \begin{equation}
    \label{eq:energy}
    \int_\domain \alpha^2 |\vv|^2 + \frac{1}{2}
    \| \nabla \vv + \nabla \vv^\intercal \|_F^2 \, \dVy
    = -\int_\bdry \left ( \ssigma(\vv) \cdot \v\nhat \right )
    \cdot \vv \, \dSy \; .
  \end{equation}
\end{lemma}

The quantity $-\ssigma(\vv)\cdot \v\nhat$ is
sometimes referred to as the surface traction.
For the modified Stokes equation, specifying
the surface traction is the natural
Neumann boundary value problem. We have
\begin{corollary}[Uniqueness \cite{biros2002embedded}]
  Let $\domain$ be a bounded (possibly multiply-connected)
  domain. There is at most one solution (up to an
  additive constant on $p$) of
  \cref{eq:app_momentum,eq:app_cty} specifying either
  the velocity field, $\vv$, or the surface traction,
  $\ssigma(\vv)\cdot \v\nhat$, on $\bdry$.
\end{corollary}

\begin{corollary}[Uniqueness, exterior \cite{biros2002embedded}]
  Let $\domain$ be the exterior of a finite collection
  of bounded, simply-connected domains. There
  is at most one solution (up to an
  additive constant on $p$) of
  \cref{eq:app_momentum,eq:app_cty} which satisfies
  $|\vv| = o(1/|\xx|)$ as $|\xx| \to \infty$, specifying either
  the velocity field, $\vv$, or the surface traction,
  $\ssigma(\vv)\cdot \v\nhat$, on $\bdry$.
\end{corollary}

While it was not necessary to the discussion in
the main text, the modified Stokes single layer
potential is useful in the discussion of jump
conditions and invertibility. 
For a given single layer density $\mmu$ defined on
$\bdry$, the potential is denoted by $\sglpot[\mmu]$
and is defined by the boundary integral
\begin{align}
  \sglpot_i[\mmu](\xx) &= \int_{\bdry} S_{ij}(\xx,\yy) \mu_j(\yy) \dSy, \label{eq:sgl_lyr_pot}
\end{align}
where $S_{ij}$ is the Stokeslet as defined in
\eqref{eq:stokeslet}. Denote the stress tensor
associated with the single layer potential
by $\ssigma_\sglpot[\mmu]$. Another layer potential of
interest corresponds to the surface traction of
the single layer potential. We denote this layer
potential by $\NN[\mmu]$, and for $\xx \in \bdry$,
it is defined by

\begin{align}
  \NN[\mmu]_i(\xx) = -\oint_\bdry T_{ijk}(\xx,\yy) \mu_j(\yy)
  \nhat_k(\xx) \dSy , \label{eg:sgl_lyr_pot_prime}
\end{align}
where the $\oint$ symbol indicates that this
integral is interpreted in the Cauchy principal
value sense. Note that this integral operator is
the transpose of the double layer operator.

The modified Stokes
single and double layer potentials satisfy jump
conditions which are analogous to the more
familiar jump conditions for harmonic layer
potentials. 

\begin{lemma}[Jump conditions]
  \label{lemma:jump-conds}
  Let $\sglpot$, $\ssigma_\sglpot$, $\NN$, and $\dblpot$ be
  the layer potentials as defined above and let $\bdry$ be a
  sufficiently smooth domain boundary
  with inward pointing normal $\v \nhat$. Then, for
  a given density $\mmu$ defined on $\bdry$,
  we have that $\sglpot \mmu$
  is continuous across $\bdry$, the exterior and interior
  limits of the surface traction of $\dblpot \mmu$ are equal,
  and for each $\xx_0 \in \bdry$,

  \begin{align}
    \lim_{h \to 0} -\ssigma_\sglpot[\mmu](\xx_0 \pm h\v\nhat(\xx_0)) \cdot \v\nhat(\xx_0)
    &= \mp \frac{1}{2} \mmu(\xx_0) + \NN[\mmu](\xx_0) \\
    \lim_{h \to 0} \dblpot[\mmu](\xx_0 \pm h\v\nhat(\xx_0)) 
    &= \pm \frac{1}{2} \mmu(\xx_0) + \dblpot[\mmu](\xx_0)    \; .
  \end{align}
  Note that the integral in the definition of $\dblpot[\mmu](\xx_0)$
  is interpreted in the Cauchy principal value sense
  when $\xx_0 \in \bdry$.
\end{lemma}

The above expressions are derived by noting that the
leading order singularity of these integral kernels
is the same as for the original Stokes case, so that
the standard jump conditions for Stokes
\cite{KimSangtae1991M:pa,Pozrikidis1992}
apply. The correspondence between the $\dblpot$ and $\NN$
operators, namely that $\dblpot = \NN^\intercal$, helps
to characterize the nullspace of $\frac{1}{2} + \dblpot$.
\begin{lemma}[Nullspace \cite{biros2002embedded}]
  \label{lemma:nullspaces}
  Suppose that $\bdry$ is the boundary of a (possibly multiply-connected)
  bounded domain and that $\v\nhat$ denotes the inward pointing
  normal. 
  Then $\dim (N(\frac{1}{2} + \dblpot)) = 1$.
  If $\xxi \in N(\frac{1}{2} + \dblpot)$ and $\xxi \ne \zzero$,
  then $\int_\bdry \xxi \cdot \v\nhat \ne 0$. 
\end{lemma}

In order to deal with the rank 1 nullspace described above,
we add the term
\begin{align}
  \mathcal W[\mmu](\xx)
  = \frac{\v\nhat(\xx)}{\int_{\bdry} \dif S} \int_{\bdry}
  \mmu(\yy) \cdot \v\nhat(\yy) \dSy \; 
\end{align}
to the double layer potential. This process is sometimes referred to
as a nullspace correction or Wielandt’s deflation \cite{biros2002embedded,Frank1958}.
We summarize this process in 
\begin{lemma}[Invertibility]
  \label{lemma:invertibility}
  Let $\mathcal{W}$ and $\dblpot$ be defined as above. Then, the
  equation
  \begin{align}
    \label{eq:dlp_append}
  \frac{1}{2}\mmu + \DD[\mmu] + \mathcal W[\mmu] = 
  \ff, \quad \xx \in \bdry 
\end{align}
is uniquely invertible. Moreover, $\mathcal{W}[\mmu]$ is
zero provided that the right hand side satisfies the
condition 

\begin{equation}
  \label{eq:compat_append}
  \int_\bdry \ff \cdot \v\nhat
  \, \dSy = 0 \; .
\end{equation}
\begin{proof}

  Suppose that $\ff$ satisfies \cref{eq:compat_append}. Taking
  the inner product of \cref{eq:dlp_append} with $\v\nhat$ and
  integrating, we obtain

  \begin{equation}
    \int_\bdry \v\nhat \cdot \left (\frac{1}{2} + \dblpot
    + \mathcal{W} \right )\mmu  = 0 \; .
  \end{equation}
  Because $\uu = \dblpot \mmu$ is divergence-free and
  $(\frac{1}{2} + \dblpot) \mmu$ gives $\uu$ on the boundary,
  the above implies that $\int \mmu \cdot \v\nhat = 0$ and
  $\mathcal{W}\mmu = \zzero$

  By the Fredholm alternative, it is sufficient to show
  that
\begin{align}
  \frac{1}{2}\mmu + \DD[\mmu] + \mathcal W[\mmu] = 
  \zzero, \quad \xx \in \bdry 
\end{align}
implies that $\mmu = \zzero$. Note that $\ff = \zzero$ certainly
satisfies \cref{eq:compat_append}. Therefore, $\int \mmu \cdot \v\nhat = 0$
and $(\frac{1}{2} + \dblpot) \mmu = \zzero$. Then, by
\Cref{lemma:nullspaces}, we have that $\mmu = \zzero$.
  
\end{proof}
\end{lemma}

\section{Free-space Green's functions of the modified Stokes equations}
\label{sec:fundamental-solutions}

We here derive the fundamental solutions to the modified Stokes equations needed for our application. Our starting point is the fundamental solution to the modified biharmonic equation,
\begin{align}
  \Delta(\Delta - \alpha^2)G(\v r) = \delta(\v r) .
  \label{eq:mod_biharm}
\end{align}
In two dimensions~\cite{Jiang2013},
\begin{align}
  G(\v r) &= -\frac{1}{2\pi\alpha^2}\left(
            \log r + K_0(\alpha r) \right),
            \label{eq:G_biharm}
\end{align}
where $r = \norm{\v r}$ and $K_n$ is a modified Bessel functions of
the second kind. We now seek the fundamental solution to the modified
Stokes equation, such that
\begin{align} 
  \alpha^2 \v u(\v r) - \Delta \v u(\v r) + \nabla p(\v r) &= \delta(\v r) \v f, 
  \label{eq:mod_stokes} \\
  \nabla \cdot \v u(\v r) &= 0 ,
                      \label{eq:divergence}
\end{align}
where $\v f$ is an arbitrary constant. Substituting $\delta(\v r)$ using \eqref{eq:mod_biharm} and taking the divergence,
\begin{align}
  \Delta p &= \nabla\cdot\Delta(\Delta - \alpha^2) G \v f ,
\end{align}
such that
\begin{align}
  p &= (\Delta - \alpha^2)\nabla G \cdot \v f .
      \label{eq:pressure_from_G}
\end{align}
Substituting \eqref{eq:mod_biharm} and \eqref{eq:pressure_from_G} into \eqref{eq:mod_stokes},
\begin{align}
  (\alpha^2-\Delta)\v u &= (\alpha^2-\Delta) \nabla (\nabla G \cdot \v f ) 
                   + \Delta(\Delta - \alpha^2)G \v f, \\
                 &= (\alpha^2-\Delta)
                          (\nabla \otimes \nabla) G \v f 
                          - (\alpha^2-\Delta) \Delta G \v f,  \\
  \v u &= \left(\nabla \otimes \nabla - \Delta I \right) G \v f, 
\end{align}
such that the stokeslet tensor $S$ is given by
\begin{align}
  S(\v r) = (\nabla\otimes\nabla -\Delta) G(\v r).
\end{align}
We now switch to index notation, following the Einstein summation
convention and using the notation $\di=\pd{}{r_i}$. The Kronecker
delta $\delta_{ij}$ is not to be confused with the Dirac delta function
$\delta(\v r)$ used above. The stokeslet is then written
\begin{align}
  S_{ij}(\v r) = \left( \di\dj - \Delta\delta_{ij} \right) G(\v r),
\end{align}
and the corresponding pressure vector is given by
\begin{align}
  \presvec_i(\v r) = \di(\Delta-\alpha^2) G(\v r) .
\end{align}
The stresslet tensor, defined as
\begin{align}
  T_{ijk}(\v r) = -\delta_{ik} \presvec_j(\v r) + \partial_k S_{ij}(\v r) + \partial_i S_{kj}(\v r),
\end{align}
is then given by
\begin{align}
  T_{ijk}= 
  -\left(\delta_{ij}\dk + \delta_{jk}\di + \delta_{ki}\dj\right) \Delta G
  + 2\di\dj\dk G 
  + \alpha^2\delta_{ik}\dj G .
\end{align}
Note that for $\alpha=0$ we recover the corresponding tensors for
regular Stokes flow, with $G$ being the fundamental solution to the
biharmonic equation.

We have that $G$ is a radial function in $r$, $G(\v r) = G(r)$. Using
that in dimension $d$
\begin{align}
  \di r &= \frac{r_i}{r}, & \di r_j &= \delta_{ij}, &
  \di r_i &= d, & r_ir_i &= r^2,
\end{align}
and
\begin{align}
  \Delta G(\v r) = G''(\v r) + (d-1)\frac{G'(r)}{r},
\end{align}
it is straightforward to derive that
\begin{align}
  \di G(\v r) &= G'(\v r) \frac{r_i}{r}, \\
  \di \dj G(\v r) &= \frac{G'(\v r)}{r} \delta_{ij} 
                    + \left( G''(\v r) 
                    - \frac{G'(\v r)}{r} \right) 
                    \frac{r_i r_j}{r^2},\\
  \begin{split}
    \di \dj \dk G(\v r) &= \left(\frac{G''(r)}{r} -
      \frac{G'(r)}{r^2}\right) \frac{ \delta_{jk} r_i + \delta_{ik}
      r_j
      + \delta_{ij} r_k }{r} \\
    & \quad + \left(G'''(r) - 3 \frac{G''(r)}{r} + 3
      \frac{G'(r)}{r^2}\right)
    \frac{r_i r_j r_k}{r^3}, 
  \end{split}\\
  \di \Delta G(\v r) &= \left(G'''(\v r) 
                       + (d-1)\frac{G''(\v r)}{r} 
                       - (d-1)\frac{G'(\v r)}{r^2} \right)
                       \frac{r_i}{r}.
\end{align}
Using these relations, we can now write down $S$, $\v\presvec$ and $T$ in
terms of $G$, $G'$ and $G''$,
\begin{align}
  S_{ij}(\v r) &= - \frac{1}{r} \left( r G''(\v r) + 
                 (d-2) G'(\v r) \right)
                 \delta_{ij}
                 + \frac{1}{r^3} \left( r G''(\v r) - G'(\v r) \right)
                 r_i r_j, \\
  \presvec_i(\v r) &= \frac{1}{r^3} \left(r^2 G'''(\v r) 
              + (d-1) r G''(\v r)
              - \left(r^2\alpha^2+d-1\right) G'(\v r) \right)
              r_i, \\
  \begin{split}
    T_{ijk}(\v r) &= \frac{1}{r^3} \left( -r^2 G'''(\v r) + (3-d) r
      G''(\v r) - (3-d) G'(\v r) \right) \left( \delta_{jk} r_i +
      \delta_{ik} r_j
      + \delta_{ij} r_k \right)\\
    & \quad + \frac{2}{r^5} \left( r^2 G'''(\v r) - 3 r G''(\v r) + 3
      G'(\v r) \right) r_i r_j r_k + \frac{\alpha^2 G'(\v r)}{r}
    \delta_{ik} r_j .
  \end{split}
\end{align}
In two dimensions (using $d=2$ and \eqref{eq:G_biharm}) we can write
\begin{align}
  S_{ij}(\v r) &= \mathcal S_1(\alpha r) \delta_{ij}
                 + \alpha^2 \mathcal S_2(\alpha r) r_i r_j,
                 \label{eq:app_stokeslet} \\
  \presvec_i(\v r) &= \frac{r_i}{2 \pi r^2}, \\
  T_{ijk}(\v r) &= \alpha^2 \mathcal T_1(\alpha r) \left( 
                  \delta_{jk} r_i 
                  + \delta_{ik} r_j 
                  + \delta_{ij} r_k \right)
                  + \alpha^4 \mathcal T_2(\alpha r) 
                  r_i r_j r_k
                  + \alpha^2 \mathcal T_3(\alpha r) \delta_{ik} r_j,
                  \label{eq:app_stresslet}
\end{align}
where
\begin{align}
  \mathcal S_1(z) &= \frac{z^2 K_0(z) + z K_1(z)-1}{2 \pi z^2},
                    \label{eq:app_stokeslet_S1} \\
  \mathcal S_2(z) &= -\frac{z^2 K_0(z)+2 z K_1(z)-2}{2 \pi  z^4
                    },
                    \label{eq:app_stokeslet_S2} \\
  \mathcal T_1(z) &= -\frac{ 2 z^2 K_0(z)+\left(z^2+4\right) z K_1(z)-4}{2 \pi  z^4},
                    \label{eq:app_stresslet_T1} \\
  \mathcal T_2(z) &= \frac{4 z^2 K_0(z)+\left(z^2+8\right) z K_1(z)-8}{\pi  z^6},
                    \label{eq:app_stresslet_T2} \\
  \mathcal T_3(z) &= \frac{z K_1(z)-1}{2 \pi  z^2}. \label{eq:app_stresslet_T3} 
\end{align}
To avoid cancellation errors, the above expressions must be evaluated
using power series for small values of the argument $z$, see
\cref{sec:power-series}.

%%%%%%%% Looks like we don't need to explicitly write the gradient in this way
%
% To evaluate the gradient of the double layer potential, we also need the kernel
% \begin{align}
%   \gradker_{ijkl} := \dl T_{ijk}= 
%   -\left(\delta_{ij}\dk + \delta_{jk}\di + \delta_{ki}\dj\right) \dl \Delta G
%   + 2\di\dj\dk\dl G 
%   + \alpha^2\delta_{ik} \dj \dl G .
% \end{align}
% To derive this we need the following additional results:
% \begin{align}
%   \begin{split}
%     &\di \dj \dk \dl G(\v r) = \left( \frac{G''(r)}{r^2} -
%       \frac{G'(r)}{r^3} \right)
%     (\delta_{ij} \delta_{kl}+\delta_{ik} \delta_{jl}+\delta_{il}\delta_{jk})  \\
%     & \quad + \left( \frac{G^{(3)}(r)}{r} - 3 \frac{G''(r)}{r^2} + 3
%       \frac{G'(r)}{r^3} \right)
%     \frac{(\delta_{kl} r_i r_j + \delta_{jl} r_i r_k+ \delta_{il} r_j r_k + \delta_{ij} r_l r_k +\delta_{jk} r_i r_l+\delta_{ik} r_j r_l)}{r^2} \\
%     & \quad + \left( G^{(4)}(r) - 6 \frac{G^{(3)}(r)}{r} + 15
%       \frac{G''(r)}{r^2} -15 \frac{G'(r)}{r^3}\right) \frac{r_i r_j
%       r_k r_l}{r^4}
%   \end{split}
% \end{align}
% \begin{align}
%       \begin{split}
%         \di \dj \Delta G(\v r) = & \left( \frac{G^{(3)}(r)}{r} + (d-1)
%           \left( \frac{G''(r)}{r^2} - \frac{G'(r)}{r^3} \right)
%         \right) \delta_{ij}
%         \\
%         &+ \left( G^{(4)}(r) + (d-2)\frac{G^{(3)}(r)}{r} + (d-1)
%           \left( -3 \frac{ G''(r)}{r^2} + 3 \frac{G'(r)}{r^3}
%           \right)\right) \frac{r_i r_j}{r^2}
%       \end{split}
% \end{align}

%%% Local Variables: 
%%% TeX-master: "nseIE"
%%% End: 
  
  \section{Power series}
\label{sec:power-series}

In the closed-form expressions for the stokeslet and the stresslet,
\cref{eq:app_stokeslet,eq:app_stresslet}, the expressions for the
functions $\mathcal S_i(z)$ and $\mathcal T_i(z)$,
\cref{eq:app_stokeslet_S1,eq:app_stokeslet_S2,eq:app_stresslet_T1,eq:app_stresslet_T2,eq:app_stresslet_T3},
are prone to cancellation errors for small arguments. To get around
this, we form their power series, by combining the coefficients of the
series expansions of the modified Bessel functions. The fundamental
expansions needed are \cite[\S10.25 and \S10.31]{NIST:DLMF}
\begin{align}
  I_0(z) &= \sum_{n=0}^\infty \frac{ \left(\frac{1}{4}z^2\right)^n}
           {(n!)^2}, \\
  I_1(z) &= \frac{z}{2}\sum_{n=0}^\infty \frac{ \left(\frac{1}{4}z^2\right)^n }{n!(n+1)!}, \\
  K_0(z) &= -\log(z/2) I_0(z) + 
           \sum_{n=0}^\infty \psi(n+1) 
           \frac{\left(\frac{1}{4}z^2\right)^n }{(n!)^2}, \\
  K_1(z) &= \frac{1}{z} + \log(z/2) I_1(z) -
           \frac{z}{4}\sum_{n=0}^\infty 
           (\psi(n+1) + \psi(n+2))
           \frac{\left(\frac{1}{4}z^2\right)^n }{n!(n+1)!},\\        
\end{align}
where $\psi$ is the digamma function, which be can evaluated recursively
starting from the Euler-Mascheroni constant $\gamma$,
\begin{align}
  \psi(1) &= -\gamma  \\
  \psi(n+1) &= \psi(n) + 1/n .
\end{align}
We write the above power series on the following compact form:
\begin{align}
  I_0(z)  &= \sum_{n=0}^\infty i_0(n) z^{2n}, &
  I_1(z)  &= \sum_{n=0}^\infty i_1(n) z^{2n+1},\\
  K_0(z)  &= \sum_{n=0}^\infty k_0(z,n) z^{2n} , &
  K_1(z)  &= \frac{1}{z} + \sum_{n=0}^\infty k_1(z,n) z^{2n+1},
\end{align}
where
\begin{align}
  i_0(n) &=  \frac{1}{4^n (n!)^2}, \\
  i_1(n) &= \frac 12 \frac{ 1 }{4^nn!(n+1)!}, \\
  k_0(z,n) &=  
    \left(\psi(n+1) -\log(z/2) \right)
    \frac{1}{4^n (n!)^2}    
    \\
  k_1(z,n) &=\left(
             2\log(z/2)
             -
             \psi(n+1) - \psi(n+2)
             \right)
             \frac{1}{4^{n+1}n!(n+1)!} 
\end{align}
In our application we only have use for the functions $\mathcal
T_i(z)$. These can now be computed for small $z$ using
\begin{align}
  \mathcal T_1(z) 
  &= -\frac{1}{2 \pi}
    \sum_{n=1}^\infty \left(
    2 k_0(z,n) + 4 k_1(z, n) + k_1(z, n-1) \right) z^{2n-2}, \\
  \mathcal T_2(z) 
  &= -\frac{1}{\pi z^4} + \frac{1}{8 \pi z^2} + \frac{1}{\pi}
    \sum_{n=2}^\infty \left(
    4 k_0(z,n) + 8 k_1(z, n) + k_1(z, n-1) \right) z^{2n-4}, \\
  \mathcal T_3(z) ,
  &= \frac{1}{2\pi}
    \sum_{n=0}^\infty k_1(z,n) z^{2n} .
\end{align}
In our code, we have empirically determined that we need to switch to
power series evaluation for $\mathcal T_i(z)$ when $z \le 1.5$, and
that it is sufficient to truncate the sums to 11 terms. For the
derivatives $\mathcal T_i'(z)$, we evaluate the derivatives of the
above series when $z \le 2$, with the sums truncated to 13 terms.

For the kernel-split quadrature described in \cref{sec:quadrature}, we
also need power series expansions of the smooth functions
$\mathcal T_i^S$ and $\mathcal T_i^L$, as defined in
\cref{sec:kernel-split}. These are straightforward to break out from
the above expansions, by splitting the $\log(z)$ factor from the
coefficients,
\begin{align}
  k_i(z, n) = k_i^S(n) + k_i^L(n) \log z .
\end{align}

%%% Local Variables: 
%%% TeX-master: "nseIE"
%%% End: 

  %\input{append-split-large-arg}
\end{appendices}

\clearpage
\bibliography{nse2d,refs-ludvig} % This is the bib file exported by Mendeley for our group

\end{document}